\documentclass[notitlepage]{article}
\usepackage{amsfonts,amssymb,amsmath,amsthm,amscd}
\usepackage{eucal}

\usepackage[all]{xy} \CompileMatrices

\usepackage[dvips]{color,graphicx}
\newtheorem{theorem}{Theorem}[subsection]
\newtheorem{proposition}[theorem]{Proposition}
\newtheorem{lemma}[theorem]{Lemma}

\theoremstyle{definition}
\newtheorem{definition}[theorem]{Definition}

\theoremstyle{remark} \newtheorem{remark}[theorem]{Remark}
\theoremstyle{definition} \newtheorem*{assumption}{Assumption}
\numberwithin{equation}{subsection}
\newcommand{\field}[1]{\ensuremath{\mathbb{#1}}} 
\newcommand{\CC}{\field{C}} \newcommand{\DD}{\field{D}}
\newcommand{\HH}{\field{H}} \newcommand{\PP}{\field{P}}
\newcommand{\RR}{\field{R}} 
\newcommand{\ZZ}{\field{Z}}
\newcommand{\complex}[1]{\mathsf{#1}} 
 \newcommand{\SSS}{\complex{S}}
\newcommand{\BBB}{\complex{B}} \newcommand{\KKK}{\complex{K}}
\newcommand{\AAA}{\complex{A}} \newcommand{\CCC}{\complex{C}}

\newcommand{\diff}[1]{\mathcal{#1}}

\newcommand{\cover}[1]{\mathcal{#1}}

\newcommand{\sheaf}[1]{\underline{\mathnormal{#1}}}

\newcommand{\cat}[1]{\mathrm{\boldsymbol{#1}}}

\newcommand{\teich}[2]{\mathcal{#1}(#2)}
\newcommand{\ttt}{\teich{T}{X}} \newcommand{\ppp}{\teich{P}{X}}
\newcommand{\sss}{\teich{S}{X}} \newcommand{\bbb}{\teich{B}{X}}
\newcommand{\gggg}{\teich{G}{X}}

 \DeclareMathOperator{\id}{id}
 \DeclareMathOperator{\Tot}{Tot}

\DeclareMathOperator{\PSL}{PSL} \DeclareMathOperator{\SO}{SO}

\newcommand{\HHH}{\mathbb{H}}

\newcommand{\del}{\partial} 
\newcommand{\delb}{\Bar\partial} 
\newcommand{\deltacheck}{\Check\delta} 
\newcommand{\delp}{\partial^\prime}
\newcommand{\delpp}{\partial^{\prime\prime}}
\newcommand{\deltapp}{\delta}

\newcommand{\var}{\boldsymbol{\delta}} 

\newcommand{\tame}[2]{\bigl(#1,#2\bigr]}
\newcommand{\abs}[1]{\lvert#1\rvert}
\newcommand{\interior}{\lrcorner}
\newcommand{\eqdef}{\overset{\mathrm{def}}{=}}

\begin{document}
\begin{titlepage}
\begin{flushright}
  SISSA 62/2000/FM
\end{flushright}
\vspace{3pc}
\begin{center}
{\LARGE
  Generating Functional in CFT on Riemann Surfaces II:
  Homological Aspects
  }\\[1.5pc]
{\large
  Ettore Aldrovandi\footnote{Current address: Department of
    Physics, Florida State University, Tallahassee, FL 32306-4350.
    \texttt{ettore@hep.fsu.edu}}\\
  S.I.S.S.A./International School for Advanced Studies\\
  Via Beirut 2/4, 34013 Trieste, Italy\\
  \texttt{ettore@fm.sissa.it}
  \\[1pc]
  Leon A. Takhtajan\\
  Department of Mathematics, SUNY at Stony Brook\\
  Stony Brook, NY 11794-3651, USA\\
  \texttt{leontak@math.sunysb.edu}
  }
\\[2pc]
\textit{Dedicated to the memory of Han Sah}
\end{center}
\vspace{2.5pc}
\begin{abstract}
  We revisit and generalize our previous algebraic construction
  of the chiral effective action for Conformal Field Theory on
  higher genus Riemann surfaces. We show that the action
  functional can be obtained by evaluating a certain Deligne
  cohomology class over the fundamental class of the underlying
  topological surface. This Deligne class is constructed by
  applying a descent procedure with respect to a \v{C}ech
  resolution of any covering map of a Riemann surface. Detailed
  calculations are presented in the two cases of an ordinary
  \v{C}ech cover, and of the universal covering map, which was
  used in our previous approach. We also establish a dictionary
  that allows to use the same formalism for different covering
  morphisms.
  
  The Deligne cohomology class we obtain depends on a point in
  the Earle-Eells fibration over the Teichm\"uller space, and on
  a smooth coboundary for the Schwarzian cocycle associated to
  the base-point Riemann surface. From it, we obtain a
  variational characterization of Hubbard's universal family of
  projective structures, showing that the locus of critical
  points for the chiral action under fiberwise variation along
  the Earle-Eells fibration is naturally identified with the
  universal projective structure.
\end{abstract}
\end{titlepage}

\tableofcontents

\section{Introduction}
\label{sec:intro}

This paper is a follow-up to our previous
paper~\cite{prev_paper}, where we presented an algebraic
construction of the chiral effective action for Conformal Field
Theory on higher genus Riemann surfaces.  The aim of the present
work is two-fold.

First, in light of the renewed interest for Classical Field
Theory~\cite{del-freed}, we present a case study for an action
functional whose construction exhibits non-trivial algebraic
properties --- the action is actually the evaluation of a certain
Deligne class. The functional is non-topological, which should be
contrasted with cases where methods of homological algebra and
algebraic topology were used to construct topological
terms~\cite{freed,gawedzki,alvarez}. Furthermore, in the recent
development of String Theory, there appear dynamical fields of a
new geometric content, such as, for example, the $B$-field.  It
is very important to find adequate geometric structures to
describe these fields and to devise suitable action
functionals~\cite{freed-wit}. Some attempts have been made at
introducing the language of \emph{gerbes} as the proper geometric
structure, at least in the lower degrees (where the language
itself makes sense). In this approach, one usually settles for a
\v Cech description relative to some open covering of the
underlying manifold. Therefore an added motivation to our work,
although we mention gerbes only in passing, was to show the
universal nature of the \v Cech paradigm for constructing action
functionals. By this we mean to develop a method which works for
general \v Cech resolutions and cohomology with respect to
\emph{arbitrary coverings}, and not just the standard open cover,
and which allows to freely change among the coverings.

This brings us to the second goal: to describe explicitly the
dependence of the chiral action functional on various default
choices, which is necessary in order to make our construction
in~\cite{prev_paper} work for arbitrary coverings. In particular,
this calls for the following.
\begin{enumerate}
\item A detailed analysis of the descent equations with respect
  to the nerve of the cover, where the use of Deligne complexes
  becomes crucial.
\item An analysis of the dependence of the chiral action on the
  choice of the projective structure on the Riemann surface.
\end{enumerate}

Recall that the choice of the universal cover for a Riemann
surface, made in~\cite{prev_paper}, yields a default choice for
the projective structure: the Fuchsian projective structure,
provided by the uniformization map. Since the universal Conformal
Ward Identity (CWI) determines the chiral action only up to a
holomorphic projective connection, the dependence of the chiral
action functional on the choice of a projective structure should
be compatible with it. Indeed, we prove this for the chiral
action ``on shell'', i.e., for solutions of the classical
equations of motion.

In order to describe the content of this paper in more detail, we
briefly recall the main results in~\cite{prev_paper}.

Let $\mu$ be a Beltrami coefficient on $\CC$ --- a smooth bounded
function $\mu$ with the property
$\Vert\mu\Vert_{\infty}=\sup_{z\in\CC} \abs{\mu(z)} < 1$ --- and
let $f$ be a solution of the Beltrami equation
\begin{equation*}
f_{\Bar z}= \mu\,f_z\,,
\end{equation*}
a self-map $f:\CC\to \CC$, unique up to  post-composition with
a M\"obius transformation. The Euclidean version of Polyakov's
action functional for two-dimensional quantum gravity~\cite{pol}
has the form
\begin{equation*}
  S[f]=2\pi i\int_{\CC} \frac{f_{zz}}{f_z}\mu_z
  \,dz\wedge d\Bar z\,,
\end{equation*}
and solves the universal Conformal Ward Identity
\begin{equation*}
  (\delb-\mu\,\del-2\,\mu_z)\frac{\var W}{\var \mu(z)}
  =\frac{c}{12 \pi} \mu_{zzz}\,,
\end{equation*}
where $W[\mu]$ is the generating functional for the vacuum chiral
conformal block, and
\begin{equation*}
  W[\mu]=-\frac{c}{96 \pi^2} S[f]\,.
\end{equation*}
Here $c$ is the central charge of the theory, and we denoted by
$\var$ the variational operator.

In~\cite{prev_paper}, we extended Polyakov's ansatz from $\CC$ to
a compact Riemann surface $X$ of genus $g > 1$, using the
following construction. Consider the universal cover $\HH\to X$,
where $\HH$ is the upper half-plane, and let $\mu$ be a Beltrami
coefficient on $\HH$, which is a pull-back of a Beltrami
coefficient on $X$ (see~\ref{sec:prel_qc} and~\cite{prev_paper},
and also~\cite{ahl,nag} for details). Depending on the extension
of $\mu$ into the lower half-plane, there exists a unique
solution $f$ to the Beltrami equation on $\HH$.  It is a map
$f:\HH \to \DD$ with the following intertwining property:
\begin{equation*}
  f\circ \Gamma = \Tilde\Gamma\circ f\,,
\end{equation*}
where $\Gamma$ is a Fuchsian group uniformizing the Riemann
surface $X$ (it is isomorphic to $\pi_1(X)$ as an abstract
group), and $\Gamma\to \Tilde\Gamma$ is an isomorphism onto a
discrete subgroup of $\PSL_2(\CC)$. The domain $\DD=f(\HH)$ is
diffeomorphic to $\HH$ and can be made equal to $\HH$ by choosing
an appropriate extension of $\mu$. In this way one gets a
deformation map $f:X\cong \Gamma\backslash \HH\to \Tilde\Gamma
\backslash \DD \cong \Tilde X$ (which is also denoted by $f$)
onto a new Riemann surface $\Tilde X$.

The de Rham complex on $\HH$ is a complex of $\Gamma$-modules for
the obvious pull-back action. The basic 2-form of Polyakov's
ansatz
\begin{equation*}
  \omega[f]=\frac{f_{zz}}{f_z}\mu_z \,dz\wedge d\Bar z
\end{equation*}
on $\HH$ is manifestly not invariant under the action of
$\Gamma$; this means that regarding $\omega[f]$ as a 0-cochain
for $\Gamma$ with values in 2-forms, its group coboundary is not
zero. Nevertheless, $\omega[f]$ can be extended to a cocycle
$\Omega[f]$ of total degree 2 living in the double complex
$\CCC^{p,q} = C^q(\Gamma,\sheaf{A}^p(\HH))$, whose total
cohomology coincides with the de Rham cohomology of $X$. Simple
integration for the genus zero case is replaced by the evaluation
over a suitable representative $\Sigma$ of the fundamental class
$[X]$ of $X$, defining
\begin{equation*}
  S[f] = \langle \Omega[f],\Sigma\rangle.
\end{equation*}
This construction~\cite{prev_paper} extends the definition of the
chiral action to a higher genus Riemann surface $X$, and the
functional $S[f]$ has the same variational properties as
Polyakov's action on the complex plane. In particular, it solves
the universal CWI, the general solution being the sum of
$W[\mu]=-c/96\pi^2\,S[f]$ and an arbitrary quadratic
differential, holomorphic with respect to the new complex
structure on $X$ determined by the Beltrami differential $\mu$.

The main advantage of working with the universal cover $\HH$ is
that one can use formulas from the genus zero case and simply
''push them onto'' the double complex $\CCC^{p,q} =
C^q(\Gamma,\sheaf{A}^p(\HH))$.\footnote{Another procedure would
  be to find a covariant version of everything on the base $X$
  (cf.~\protect\cite{laz,zucchini}), but this introduces
  additional ``background'' structures with no direct bearing to
  the complex and algebro-topological structures of $X$.}
However, working with the universal cover uses several default
choices, as follows.
\begin{itemize}
\item The groups $\Gamma$ and $\Tilde\Gamma$ are discrete
  subgroups of $\PSL_2(\RR)$ and $\PSL_2(\CC)$ respectively, so
  that local sections to the covering maps $\HH\to X$ and $\DD\to
  \Tilde X$ are projective structures subordinated to the complex
  structures of $X$ and $\Tilde X$, respectively. These
  projective structures are inherent in the choice of $\HH$ as a
  cover, and they do not appear explicitly in the expression for
  the total cocycle $\Omega[f]$.
\item $H^3(X,\CC)=0$ has to be invoked to close the descent
  equations leading from $\omega$ to the total cocycle $\Omega$.
  This fact can be interpreted as the vanishing of an obstruction
  or, in other words, as an integrability property for the
  problem of choosing integration constants to the last descent
  equation. An element of arbitrariness is introduced in the
  explicit computation of $\Omega$ by choosing a shift of a
  $\CC$-valued 3-cochain in this equation to turn it into \v Cech
  coboundary.
\item A specific choice of logarithm branches was made
  in~\cite{prev_paper}.
\end{itemize}

The analysis of this construction shows that what we have used
were not some specific features of the universal cover $\HH\to
X$, but rather its algebraic properties relative to the double
complex $\CCC^{p,q}$: the facts that $\HH$ is contractible, and
that $\Gamma$ is cohomologically trivial with respect to modules
of smooth forms on $\HH$. These are precisely the properties of a
``good'' cover~\cite{bott-tu}, one for which the \v Cech-de Rham
double complex computes cohomology groups for both theories.

As in~\cite{prev_paper}, start with the deformation map
$f:X\to\Tilde X$, defined, say, as the solution of the Beltrami
equation on $X$.  It is natural to ask whether it is possible to
carry out the same scheme as with $\HH$ with respect to a
different cover of $X$, for example an ordinary open cover
$\cover{U}_X=\{U_i\}_{i\in I}$ of $X$, with the requirement that
it should allow for a change of covering morphism without
changing the formalism. This is achieved by considering, for a
given covering map $U\to X$ and a sheaf $\sheaf{F}$, or complex
of sheaves $\sheaf{F}^\bullet$ on $X$, its \v Cech cohomology
$\Check{H}^\bullet(U\to X;\sheaf{F})$, or hypercohomology
$\Check{\HH}^\bullet (U\to X;\sheaf{F}^\bullet)$, respectively.
The framework of the universal cover is retrieved from the
observation that group cohomology for $\Gamma$ is \v Cech
cohomology for the covering $\HH\to X$.

Our main difference from~\cite{prev_paper} is the use of the
Deligne complex instead of the simpler de Rham complex. In
particular, introducing the smooth de Rham sheaves
$\sheaf{A}_X^\bullet$, we work with the Deligne complex of length
3: $\ZZ(3)^\bullet_\mathcal{D}: \ZZ(3)\overset{\imath}{\to}
\sheaf{A}^0_X \overset{d}{\to} \sheaf{A}^1_X\overset{d}{\to}
\sheaf{A}^2_X$, where $\ZZ(3) \eqdef (2\pi i)^3\ZZ$, and apply
the same procedure as before. Namely, we form the double complex
$\CCC^{p,q}=\Check{C}^q(U\to X;\ZZ(3)^p_\mathcal{D})$, localize
the Polyakov's 2-form $\omega$ to $U$ as an element of degree
$(3,0)$ in this complex\footnote{There is a degree shift caused
  by the insertion of the integers at degree zero in the Deligne
  complex.}, and perform the usual descent calculations. The
latter procedure was first introduced into mathematical physics
in~\cite{FS}. Specifically, we solve for elements $\theta$ and
$\Theta$ of degree $(2,1)$ and $(1,2)$, respectively, satisfying
equations $\deltacheck\omega = d\theta$ and $\deltacheck\theta =
d\Theta$, with $\deltacheck\Theta\in \ZZ(3)$, where $\deltacheck$
is the \v Cech coboundary operator.  It is crucial that these
equations are solvable due to the vanishing of the tame symbol
$\tame{TX}{TX}$ in \emph{holomorphic} Deligne cohomology. As a
result, starting from Polyakov's 2-form $\omega[f]$ we obtain a
cocycle $\Omega[f]$ of total degree 3 in the total complex
$\Tot\CCC^{\bullet,\bullet}$. This constitutes the first result
of the paper, Proposition~\ref{prop:tot_cocycle}. Note that it is
convenient, for a regular open cover $\cover{U}_X$, to consider
the most general form of the bulk term for the Polyakov's action,
given by adding a smooth projective connection $h$ to the local
basic 2-form for genus 0:
\begin{equation*}
  \omega[f]=\frac{f_{zz}}{f_z}\mu_z \,dz\wedge d\Bar z +2\,\mu\,
  h \,dz\wedge d\Bar z\,.
\end{equation*}
Here $z$ is a local coordinate for $U \in \cover{U}_X$, and $h$ a
representative in $U$ of a smooth projective connection on $X$
--- a smooth coboundary for the usual Schwarzian cocycle relative
to the cover $\cover{U}_X$. The space $\mathcal{Q}(X)$ of all
such coboundaries is an affine space over the vector space of
smooth quadratic differentials on $X$. On $\HH$, the pull-back of
a projective connection is a quadratic differential. See sections
\ref{sec:prel_sheaves}-\ref{sec:prel_coverings} and
\ref{sec:action_cech}-\ref{sec:local_cocycle} for details.

In section~\ref{sec:dictionary}, we translated the generalized
\v{C}ech formalism for the universal cover $\HH\to X$ into group
cohomology for $\Gamma\cong\pi_1(X)$, so that
Proposition~\ref{prop:tot_cocycle} translates into
Proposition~\ref{prop:transl_tot_cocycle}, thus refining the
corresponding results in~\cite{prev_paper}.

For the construction of the action functional, we need to
evaluate the cocycle $\Omega[f]$ against the fundamental class
$[X]$ of a Riemann surface $X$, which we represent as a cycle
$\Sigma$ in a homological double complex $\SSS_{p,q}=S_p(N_qU)$
of singular $p$-simplices in the $q+1$-fold product of $U$ with
itself. Using the pairing $\langle ~, ~\rangle$ between Deligne
cocycles and cycles, which is well-defined because $\dim
X=2=3-1$, we can define
\begin{displaymath}
  S[f]=\langle\Omega[f],{}^\backprime\Sigma\rangle\,,
\end{displaymath}
where ${}^\backprime\Sigma$ is the shift of the cycle $\Sigma$ so
that it has total homological degree 3. Due to the insertion of
integers into the Deligne complex, the pairing $\langle
~,~\rangle$ is well defined only modulo $\ZZ(3)$, so that the
action functional $S[f]$ is well-defined only modulo $\ZZ(3)$.
Using the exponential map $z\to \exp\{z/(2\pi i)^2\}$, that
identifies $\CC/\ZZ(3)$ with $\CC^*$, one can replace the complex
$\ZZ(3)^\bullet_\mathcal{D}$ with $\sheaf{A}^*_X
\xrightarrow{d\log} \sheaf{A}^1_X\overset{d}{\to} \sheaf{A}^2_X$
and resets all degrees by one, so that cocycle $\Omega$ would
correspond to a cocycle $\Psi$ of total degree 2. The
corresponding pairing $\langle ~,~\rangle_m$ will be now
multiplicative and single-valued, with values in $\CC^*$. As a
result, the single-valued functional
\begin{displaymath}
  A[f]=\langle \Psi[f],\Sigma\rangle_m
  =\exp\{S[f]/(2 \pi i)^2\}
\end{displaymath}
is the exponential of the action, which is quite natural since we
are dealing with an effective action in QFT. Details of this
construction are presented in sections~\ref{sec:prel_sheaves}
and~\ref{sec:prel_fundamental}.

In section~\ref{sec:properties} we prove the independence of the
the functional $A[f]$ from the choices of logarithm branches,
establish its relations with Bloch dilogarithms, and show that it
can be considered as $\CC^*$-torsor.

{}The second result of the paper should be understood from the
view-point of Classical Field Theory. Let $\bbb\to\ttt$ be the
Earle-Eells principal fibration over the Techm\"uller space
$\ttt$.  The total space $\bbb$ of this fibration is the unit
ball in the $L^{\infty}$ norm in the space of all smooth Beltrami
differentials on $X$. To every $\mu\in\bbb$ there corresponds a
deformation map $f(\mu): X\to\Tilde X$, a solution of the
Beltrami equation on $X$, uniquely determined by the condition
that when pulled back to the universal cover $\HH$, it gives a
Fuchsian deformation, i.e.~$f(\HH)=\HH$. This allows to consider
the functional $A[f]$ as a map
$A:\mathcal{Q}(X)\times\bbb\to\CC^*$.

When studying the variational problem for the functionals $S[f]$
and $A[f]$, we consider the deformation map $f$ as the dynamical
field, and the projective connection $h$ as an external field,
with the problem to compute the variation with respect to $f$.
Geometrically, these variations are tangent vectors to $\bbb$,
and are of two types, depending on whether they deform the
complex structure of $X$ or not, i.e., whether the associated
Kodaira-Spencer cocycle (see section~\ref{sec:variation}) is
holomorphically trivial or not. In the former case, the
variations correspond to vertical tangent vectors to the
Earle-Eells fibration $\bbb \to \ttt$, and here we consider only
these variations.

One needs to show that this variational problem is well-defined
even though the action itself is not expressed in terms of a
simple integration over $X$ of a 2-form. In ``physical''
terminology, the bulk term given by the 2-form $\omega$ is a
multi-valued one, and we prove in Theorem~\ref{thm:variation}
that the variation of the action depends solely on the variation
of the bulk term and is a well-defined 2-form on $X$. We give two
proofs of this result. The first one is based on a careful
analysis of the descent equations for the variations of all
components of the Deligne cocycle $\Omega[f]$. The second proof,
albeit in a sketchy form, shows that this result is, in fact,
more general, and depends only on descent properties of the
variational bicomplex.  Takens'
results~\cite{takens,del-freed,zuckerman} are essential in this
context. We plan to return to this result with more details in a
more general situation, not limited to dimension 2, elsewhere.

However, this result holds only thanks to the good gluing
properties of the variations, which follow from the triviality of
the Kodaira-Spencer cocycle, and this formalism can not be
directly applied to the case of general variations. In this
respect, we point out that there was an error in the computation
of general variation in the universal cover
formalism~\cite{prev_paper}. While a brute-force calculation
would achieve the goal, we prefer to defer it until the
development of the proper treatment of the variational formalism
for multi-valued actions, where variational bicomplex(es) glue in
a more complicated way due to the non-vanishing of the
deformation class.

Returning to the present paper, we also give a geometric
interpretation of Theorem~\ref{thm:variation}. It states that at
critical points under vertical variations of the dynamical field
$f$, the external field --- the smooth projective connection $h$
--- is holomorphic with respect to the complex structure on $X$
defined by the deformation map $f$. In section
\ref{sec:projective}, we reformulate this by saying that the
space of critical points coincides with the pull-back to $\bbb$
of Hubbard's universal projective structure $\ppp\to\ttt$,
studied in~\cite{hubbard,nag}.

The paper is organized as follows. In
section~\ref{sec:preliminaries} we set up some necessary tools.
In particular, we give a brief tour of Deligne complexes and
explain the \v Cech formalism with respect to a covering $U\to
X$. We also present the minimum amount of formulas necessary to
perform the evaluation over representatives of the fundamental
class $[X]$. A more in-depth presentation would have led us
through a rather long detour from the main line of the paper,
therefore we provide it in the appendix,
in~\ref{sec:app_fund_class}.  Sections~\ref{sec:Action}
and~\ref{sec:variation_proj} comprise the main body of the paper.
After some general remarks in~\ref{sec:action_cech}
and~\ref{sec:local_cocycle}, we construct the representative
cocycle $\Omega[f]$, using \v Cech formalism with respect to an
open cover. We analyze the changes under redefinition of the
logarithm branches and of the trivializing coboundary for the
tame symbol $\tame{T_X}{T_X}$ in~\ref{sec:log_dep}.
In~\ref{sec:dictionary}, we present our construction in the form
suitable for coverings $U\to X$ other than the open one
$\cover{U}_X$, and in particular translate everything in terms of
$U=\HH$.  Finally, in~\ref{sec:variation} we discuss the
variation of the action. After a brief reminder of some basic
notions about families of projective structures
in~\ref{sec:projective}, we present in~\ref{sec:vert_variation} a
geometric interpretation of the vertical variation of the action.

\section{Preliminaries and notations}
\label{sec:preliminaries}

\subsection{Quasi-conformal maps and deformations}
\label{sec:prel_qc}

Let $X$ be compact Riemann surface of genus $g>1$. A Riemann
surface is called marked, if a system of standard generators of
its fundamental group $\pi_1(X)$ is chosen (up to an inner
automorphism).  Let $\ttt$ be the Teichm\"uller space of marked
compact Riemann surfaces of genus $g$, with base point the
Riemann surface $X$. It is defined as the set of equivalence
classes of orientation preserving diffeomorphisms
\begin{displaymath}
  f:X\longrightarrow \Tilde X,
\end{displaymath}
where the triples $[X,f_1,\Tilde{X}_1]$ and $[X,f_2,\Tilde{X}_2]$
are said to be equivalent if the map $f_2\circ f^{-1}_1$ is
homotopic to a conformal mapping of $\Tilde{X}_1$ onto
$\Tilde{X}_2$. It is well-known (see, e.g.,~\cite{nag}), that
$\ttt$ is a smooth manifold of real dimension $6g-6$, and it
admits a complex structure.

For any quasi-conformal map $f:X\rightarrow\Tilde X$, let
$\mu=\mu(f)$ be the Beltrami differential for $X$ associated to
$f$. It is a section of $TX\otimes {\Bar T}X^*$, where $TX$ is
the \emph{holomorphic} tangent bundle of $X$, satisfying the
Beltrami equation
\begin{displaymath}
  \delb f = \mu\, \del f\, ,
\end{displaymath}
where $\del=\del/\del z, \delb=\del/\del \bar z$.  Conversely, if
a $C^\infty$ Beltrami differential $\mu$ has $L^\infty$-norm less
than one, $\Vert\mu\Vert_{\infty}<1$, then the Beltrami equation
is solvable and its solution $f$ is a diffeomorphism.

Denote by $A^{-1,1}(X)=\Gamma (X,TX\otimes {\Bar TX}^*)$ the
vector space of all smooth Beltrami differentials for $X$, and by
$\bbb$ the open unit ball in $A^{-1,1}(X)$ with respect to the
$L^\infty$-norm. It is known that $\bbb$ is the total space of a
smooth infinite-dimensional principal fibration over $\ttt$ with
structure group $\gggg$, the group of all orientation preserving
diffeomorphisms of $X$ isotopic to the identity~\cite{ee,nag}.
Briefly, for every $\mu\in\bbb$ we lift it to the universal cover
$\HH$ and consider the solution $f(\mu)$ of the Beltrami equation
on $\HH$ with the condition that $f(\HH)=\HH$. Such an $f$ exists
and is unique up to a post-composition with M\"obius automorphism
of $\HH$. If $g\in\gggg$, then $\mu^g\eqdef \mu(f\circ g)$.

This provides an identification between the description of the
Teichm\"uller space as the space of equivalence classes of the
triples $[X,f,\Tilde X]$ with fixed $X$, and as the quotient of
$\bbb$ by $\gggg$.

For any $\mu\in \bbb$ denote by $[\mu]$ the corresponding element
in $\ttt$ and by $f(\mu):X\to X_\mu$ the resulting deformation of
$X$.  Though actually $X_\mu$ depends only on the class $[\mu]$,
we suppress this in the notation, and whenever the element $\mu$
is fixed, or clear from the context, we denote $X_\mu$ by $\Tilde
X$, as above.

Let $A^{p,q}(X)=\Gamma(X,{TX^*}^{\otimes p}\otimes {{\Bar
    TX}^*}^{\otimes q})$ be the space of $C^\infty$ tensors of
weight $(p,q)$, with the proviso that we take the tangent bundle
whenever either $p$ or $q$ is negative (like $A^{-1,1}(X)$ for
Beltrami differentials).  Denote by $\sheaf{A}^{p,q}_X$ the
corresponding sheaves of sections. It is well-known that the
operator
\begin{equation}
  \label{eq:delbmu}
  \delb_\mu = \delb -\mu\del -k\del\mu : \sheaf{A}^{k,0}_X\to
  \sheaf{A}^{k,1}_X
\end{equation}
is the $\delb$-operator for the complex structure determined by
$\mu$ --- the pull-back by $f$ of the complex structure on
$X_\mu$. This gives rise to the exact sequence
\begin{equation*}
  0 \to A^{-1,0}(X) \xrightarrow{\delb_\mu} A^{-1,1}(X) \to
  H^1(X_\mu,\Theta_\mu) \to 0\,,
\end{equation*}
where $\Theta_\mu$ is the tangent sheaf of $X_\mu$, which is
isomorphic to
\begin{equation*}
  0 \to T_{\mu}(\bbb /\ttt) \to T_{\mu}(\bbb) \to T_{[\mu]}(\ttt)
  \to 0\, ,
\end{equation*}
and provides the canonical identification
$T_{[\mu]}(\ttt)=H^1(X_\mu,\Theta_\mu)$ (see, e.g.~\cite{nag}).

\subsection{Sheaves and Deligne complexes}
\label{sec:prel_sheaves}

For any smooth manifold $M$, we denote by $\sheaf{A}^p_M$ the
sheaf of smooth complex-valued $p$-forms on $M$, and by $A^p(M)$
the corresponding spaces of global sections. Then
$\sheaf{A}^0_M\equiv \sheaf{A}_M$, the sheaf of smooth
complex-valued functions. When $M$ is complex, we denote by
$\sheaf{\Omega}^p_M$ the sheaves of holomorphic $p$-forms. In
particular, $\sheaf{\Omega}^0_M\equiv\mathcal{O}_M$, the
structure sheaf.

Recall that the \emph{hypercohomology} groups
$\HHH^p(M,\sheaf{F}^\bullet)$ of a complex of sheaves
\begin{displaymath}
  \sheaf{F}^\bullet:\sheaf{F}^0 \longrightarrow \sheaf{F}^1
  \longrightarrow \cdots
\end{displaymath}
are defined as the cohomology groups of the total complex of a
suitable resolution $\sheaf{I}^{\bullet,\bullet}$ of the complex
$\sheaf{F}^\bullet$. In practice, one usually takes a \v{C}ech
resolution relative to some (sufficiently fine) cover
$\cover{U}_M$ of $M$ and considers the double complex
\begin{displaymath}
  \CCC^{p,q}\eqdef \Check{C}^q(\cover{U}_M,\sheaf{F}^p)\,.
\end{displaymath}
The hypercohomology $\HHH^p(M,\sheaf{F}^\bullet)$ is computed by
taking $H^p(\Tot\CCC^{\bullet,\bullet})$, with the convention
that the total differential $D$ in degree $(p,q)$ is given by
$D=d + (-1)^p \deltacheck$, where $d$ is the differential in the
complex $\sheaf{F}^\bullet$ and $\deltacheck$ is the differential
in the \v{C}ech direction. Furthermore, two complexes
$\sheaf{F}^\bullet$ and $\sheaf{G}^\bullet$ are said to be
\emph{quasi-isomorphic} if there is a morphism $\sheaf{F}^\bullet
\to \sheaf{G}^\bullet$ inducing an isomorphism of their
cohomology sheaves: $H^\bullet(\sheaf{F})\xrightarrow{\sim}
H^\bullet(\sheaf{G})$. The standard (spectral sequence) argument
implies that their hypercohomology groups are the same. We will
apply this machinery to the case when the complex
$\sheaf{F}^\bullet$ is a smooth Deligne complex.

The use of Deligne complexes is nowadays fairly common, so we
just recall the notations and a few basic facts needed in the
sequel. It is convenient to use the ``algebraic geometers'
twist'' and set $\ZZ(p) \eqdef (2\pi i)^p\ZZ$.
Following~\cite{esn-vie,bry_book} we have:
\begin{definition}
  Let $M$ be a smooth manifold. The following complex of sheaves
  \begin{equation*}
    \label{eq:deligne_complex}
    \ZZ(p)^\bullet_\mathcal{D}: \ZZ(p)_M
    \overset{\imath}{\longrightarrow} \sheaf{A}_M
    \overset{d}{\longrightarrow} \sheaf{A}^1_M
    \overset{d}{\longrightarrow} \dots
    \overset{d}{\longrightarrow} \sheaf{A}^{p-1}_M
  \end{equation*}
  is called the \emph{smooth Deligne complex}. The \emph{smooth
    Deligne cohomology groups} of $M$ --- denoted by
  $H^q_\mathcal{D}(M,\ZZ(p))$ --- are the hypercohomology groups
  $\HHH^q(M,\ZZ(p)^\bullet_\mathcal{D})$.
\end{definition}
\begin{remark}
  $\ZZ(p)$ is placed in degree zero and the degree of each term
  $\sheaf{A}^r_M$ in $\ZZ(p)^\bullet_\mathcal{D}$ is $r+1$. The
  first differential is just the inclusion $\imath$ of $\ZZ(p)$
  in $\sheaf{A}_X$, while $d$ is the usual de Rham differential.
  The complex is truncated to zero after degree $p$. An
  equivalent definition of the Deligne complex is presented in
  the appendix, cf.~\ref{sec:cones}.
\end{remark}
The exponential map $e:\sheaf{A}_M\to \sheaf{A}^*_M$, $f\mapsto
\exp(f/(2\pi i)^{p-1})$, induces a quasi-isomorphism
\begin{displaymath}
  \ZZ(p)^\bullet_\mathcal{D}\cong (\sheaf{A}^*_M
  \xrightarrow{d\log}\sheaf{A}^1_M \overset{d} {\longrightarrow}
  \dots \overset{d}{\longrightarrow} \sheaf{A}^{p-1}_M)[-1]\,,
\end{displaymath}
where $[-1]$ denotes the operation of shifting a complex one step
to the right. Namely, for a complex $\sheaf{F}^\bullet$ the
complex $\sheaf{F}^\bullet[-1]$ is defined as
$\sheaf{F}[-1]^k=\sheaf{F}^{k-1}$, with $d_{[-1]}=-d$.

To prove this quasi-isomorphism, observe that the non zero
cohomology sheaves of the complex $\ZZ(p)^\bullet_\mathcal{D}$
are $\CC_M/\ZZ(p)_M$ and $\sheaf{A}^{p-1}_M/d\sheaf{A}^{p-2}_M$,
located in degree $1$ and $p$, respectively. Next, consider the
standard exponential exact sequence
$0\longrightarrow\ZZ(p)_M\overset{i}{\longrightarrow} \sheaf{A}_M
\overset{e}{\longrightarrow}\sheaf{A}^*_M\longrightarrow1$,
implying the following commutative diagram
\begin{equation*}
  \begin{CD}
    \ZZ(p)_M @>\iota>> \sheaf{A}_M @>d>> \sheaf{A}^1_M @>d>>
    \cdots @>d>> \sheaf{A}^{p-1}_M\\ && @VeVV @VVV & & @VVV\\ &&
    \sheaf{A}^*_M @>-d\log>> \sheaf{A}^1_M @>-d>> \cdots @>-d>>
    \sheaf{A}^{p-1}_M
  \end{CD}
\end{equation*}
where the first vertical arrow on the left is the exponential
map, and the others are given by multiplication by
$(-1)^{k-1}/(2\pi i)^{p-1}$ in degree $k$. Now it is obvious that
the two complexes have the same cohomology sheaves (by
identifying $\CC/\ZZ(p)\cong \CC^*$ through the exponential map)
and therefore have the same hypercohomology groups, up to an
index shift: $H^q_\mathcal{D}(M,\ZZ(p))\cong
\HHH^{q-1}(M,\sheaf{A}^*_M\to \sheaf{A}^1_M\to
\cdots\to\sheaf{A}^{p-1}_M)$.

\begin{remark} \label{rem:del_top_dim}
  In general, the truncation of the Deligne complex
  $\ZZ(p)^\bullet_\mathcal{D}$ after degree $p$ is fundamental.
  However, when $\dim M=p-1$, this truncation is irrelevant. In
  other words, when the length of the complex coincides with the
  dimension, $\ZZ(p)^\bullet_\mathcal{D}$ becomes an augmented de
  Rham complex: $\ZZ(p)_M\to\sheaf{A}^\bullet_M$~\cite{esnault}.
  Therefore the only non trivial cohomology sheaf occurs in
  degree 1, and $\ZZ(p)^\bullet_\mathcal{D}$ becomes
  quasi-isomorphic to $\CC_M/\ZZ(p)_M[-1]$. As a result,
  \begin{equation*}
    H^q_\mathcal{D}(M,\ZZ(p))\cong H^{q-1}(M,\CC/\ZZ(p))\cong
    H^{q-1}(M,\CC^*) \,,
  \end{equation*}
  where the latter isomorphism is given by the exponential map.
\end{remark}

Working out explicitly the first cohomology groups, one gets the
following isomorphisms: $H^1_\mathcal{D}(M,\ZZ(1)) \cong
H^0(M,\sheaf{A}^*_M)$ --- the multiplicative group of global
invertible functions --- $H^2_\mathcal{D}(M,\ZZ(1)) \cong
H^1(M,\sheaf{A}^*_M)$ --- the group of isomorphism classes of
smooth line bundles --- and $H^2_\mathcal{D}(M,\ZZ(2)) \cong
\HHH^1(M,\sheaf{A}^*_M \rightarrow \sheaf{A}^1_M)$ --- the group
of isomorphism classes of line bundles with connection. Higher
Deligne cohomology groups describe more complicated higher
geometric structures --- e.g., \emph{gerbes} and
$2$-\emph{gerbes}.

When $M$ is complex, there is an entirely analogous definition
for the \emph{holomorphic Deligne complex:}
\begin{equation*}
  \ZZ(p)^\bullet_{\mathcal{D},\mathit{hol}}: \ZZ(p)_M
  \overset{\imath}{\longrightarrow} \sheaf{\Omega}_M
  \overset{d}{\longrightarrow} \sheaf{\Omega}^1_M
  \overset{d}{\longrightarrow} \dots \overset{d}{\longrightarrow}
  \sheaf{\Omega}^{p-1}_M,
\end{equation*}
with the \emph{holomorphic Deligne cohomology groups}
$H^\bullet_{\mathcal{D},\mathit{hol}}(M,\ZZ(p))$ being the
hypercohomology groups of the complex
$\ZZ(p)^\bullet_{\mathcal{D},\mathit{hol}}$.

Many of the formal properties of the smooth Deligne complex are
also valid in the holomorphic category. In particular, there is
the exponential quasi-isomorphism
\begin{displaymath}
  \ZZ(p)^\bullet_{\mathcal{D},\mathit{hol}}\cong
  (\sheaf{\Omega}^*_M \xrightarrow{d\log}\sheaf{\Omega}^1_M
  \overset{d} {\longrightarrow} \dots
  \overset{d}{\longrightarrow} \sheaf{\Omega}^{p-1}_M)[-1]\,,
\end{displaymath}
since non trivial cohomology sheaves of these complexes occur
only in degrees $1$ and $p$ and coincide, which implies the
isomorphism in the hypercohomology, so that
$H^q_{\mathcal{D},\mathit{hol}}(M,\ZZ(p))\cong
\HHH^{q-1}(M,\sheaf{\Omega}^*_M\to \sheaf{\Omega}^1_M\to
\cdots\to \sheaf{\Omega}^{p-1}_M)$. When $\dim_\CC M=p-1$ the
truncation becomes irrelevant and
$\ZZ(p)^\bullet_{\mathcal{D},\mathit{hol}}$ is just $\ZZ(p)_M\to
\sheaf{\Omega}^\bullet_M$. Therefore, thanks to the exactness of
the holomorphic de Rham complex,
$\ZZ(p)^\bullet_{\mathcal{D},\mathit{hol}}$ is also
quasi-isomorphic to $\CC_M/\ZZ(p)_M[-1]$, and we have
\begin{displaymath}
  \HH^{q}(M,\ZZ(p)^\bullet_{\mathcal{D},\mathit{hol}}) \cong
  H^{q-1}(M,\CC/\ZZ(p))\cong H^{q-1}(M,\CC^*).
\end{displaymath}
In particular, when $M$ is a Riemann surface $X$ and $p=2$ we
have, for obvious dimensional reasons
\begin{displaymath}
  \HH^3(X,\ZZ(2)^\bullet_{\mathcal{D},\mathit{hol}}) \cong
  H^2(X,\CC^*)\cong\CC^*
\end{displaymath}
and
\begin{displaymath}
  \HH^4(X,\ZZ(2)^\bullet_{\mathcal{D},\mathit{hol}}) \cong
  H^3(X,\CC^*)=0.
\end{displaymath}
These elementary facts will play a major role in the
constructions in sect.~\ref{sec:Action}

There is a cup product $\cup: \ZZ(p)^\bullet_\mathcal{D} \otimes
\ZZ(q)^\bullet_\mathcal{D} \rightarrow
\ZZ(p+q)^\bullet_\mathcal{D}$ given by~\cite{esn-vie,bry_book}:
\begin{displaymath}
  f\cup g =
  \begin{cases}
    f\cdot g & \deg f = 0\,,\\ f\wedge dg & \deg f\geq 0\quad
    \text{and}\quad \deg g =q\,,\\ 0 &\text{otherwise,}
  \end{cases}
\end{displaymath}
and induced product in cohomology: $\cup :
H^r_\mathcal{D}(M,\ZZ(p)) \otimes H^s_\mathcal{D}(M,\ZZ(q)) \to
H^{r+s}_\mathcal{D}(M,\ZZ(p+q))$. Note that since Deligne
cohomology is defined using resolutions of complexes of sheaves,
one has to take into account the appropriate sign rules.  That
is, for two complexes $\sheaf{F}^\bullet$ and $\sheaf{G}^\bullet$
one forms the double complexes
\begin{displaymath}
  \CCC^{p,q}(\sheaf{F}) =
  \Check{C}^q(\cover{U}_X,\sheaf{F}^p)\quad \text{and} \quad
  \CCC^{r,s}(\sheaf{G}) = \Check{C}^s(\cover{U}_X,\sheaf{G}^r)
\end{displaymath}
and defines the cup product
\begin{displaymath}
  \cup : \CCC^{p,q}(\sheaf{F}) \otimes \CCC^{r,s}(\sheaf{G})
  \longrightarrow \Check{C}^{q+s}(\cover{U}_X,\sheaf{F}^p\otimes
  \sheaf{G}^r) \subset \CCC^{p+r,q+s}(\sheaf{F}\otimes\sheaf{G})
\end{displaymath}
of two elements $\{f_{i_0,\dots,i_q}\}\in \CCC^{p,q}(\sheaf{F})$
and $\{g_{j_0,\dots,j_s}\} \in \CCC^{r,s}(\sheaf{G})$ by
\begin{displaymath}
  (-1)^{qr}\,f_{i_0,\dots,i_q}\otimes
  g_{i_q,i_{q+1},\dots,i_{q+s}} \,.
\end{displaymath}
In this formula, one could replace the $\otimes$ by any other
product $\sheaf{F}^\bullet \otimes \sheaf{G}^\bullet \rightarrow
(\sheaf{F}^\bullet \cup \sheaf{G}^\bullet)$, in particular by the
cup product for Deligne complexes, introduced above.

Brylinski and McLaughlin~\cite{bry_mcl2} spell out several cup
products for the first few degrees representing interesting
symbol maps. We will use one of them later, so here we recall its
construction.

As already observed, $H^2_\mathcal{D}(M,\ZZ(1))$ corresponds to
the group of smooth line bundles on $M$. Working out details of
the \v{C}ech resolution relative to the \v{C}ech cover
$\cover{U}_M=\{U_i\}_{i\in I}$, one finds that a class in
$H^2_\mathcal{D}(M,\ZZ(1))$ is represented by the cocycle
$(f_{ij},m_{ijk})$, where $f_{ij} \in \Gamma(U_i\cap
U_j,\sheaf{A}^0_M)$ and $m_{ijk}\in \Gamma(U_i\cap U_j\cap
U_k,\ZZ(1)_M)$ are subject to the relations:
\begin{gather*}
  f_{jk} - f_{ik} + f_{ij}= m_{ijk}\,,\\ m_{jkl} - m_{ikl} -
  m_{ijl} + m_{ijk}= 0\,.
\end{gather*}
Thus $g_{ij}=\exp f_{ij}$ is a \v{C}ech 1-cocycle with values in
invertible functions, as expected.

Consider now two line bundles $L$ and $L^\prime$ over $M$,
represented by cocycles $(f_{ij},m_{ijk})$ and
$(f^\prime_{ij},m^\prime_{ijk})$, respectively. Their cup
product, to be denoted by the ``tame'' symbol
$\tame{L}{L^\prime}$ (see, e.g.,~\cite{bry_mcl2}), is an element
of $H^4_\mathcal{D}(M,\ZZ(2))$, represented by the cocycle
\begin{displaymath}
    (-f_{ij}\,df^\prime_{jk}\,,\,m_{ijk}f^\prime_{kl}\,,\,
    m_{ijk}m^\prime_{klp})\,.
\end{displaymath}

A similar interpretation holds for the holomorphic Deligne
cohomology.  In particular,
$H^2_\mathcal{D,\mathit{hol}}(M,\ZZ(1))$ corresponds to the group
of holomorphic line bundles on $M$, and the cup product of two
such line bundles is $\tame{L}{L^\prime}\in
H^4_\mathcal{D,\mathit{hol}}(M,\ZZ(2))$.  When $\dim_{\CC}M=1$,
according to the previous remark, the cup product of two
holomorphic line bundles is a trivial cocycle:
$\tame{L}{L^\prime}=0$.

In this paper our main emphasis will be on smooth Deligne
cohomology in degree three. With respect to the \v{C}ech cover
$\cover{U}_M$, a class in $H^3_\mathcal{D}(M,\ZZ(3))$ is
represented by the total cocycle
$(\omega_i,a_{ij},f_{ijk},m_{ijkl})$, where $\omega_i \in
\Gamma(U_i,\sheaf{A}^2_M)$, $a_{ij} \in\Gamma (U_i\cap
U_j,\sheaf{A}^1_M)$, $f_{ijk}\in\Gamma(U_i\cap U_j\cap U_k,
\sheaf{A}^0_M)$, and $m_{ijkl} \in \Gamma(U_i\cap U_j\cap U_k\cap
U_l,\ZZ(3)_M)$ are subject to the relations:
\begin{equation}
  \label{eq:del_3_cocycle}
  \begin{aligned}
    \omega_j-\omega_i=d\,a_{ij}, &\quad a_{jk} -a_{ik} +a_{ij}=
    -d f_{ijk},\\ \deltacheck f_{ijkl}=m_{ijkl}, &\quad
    \deltacheck m_{ijklp} = 0\,.
  \end{aligned}
\end{equation}
According to~\cite{bry_book,bry_mcl1},
$H^3_\mathcal{D}(M,\ZZ(3))$ is the group of isomorphism classes
of \emph{gerbes} on $M$, equipped with \emph{connective
structure} described by $\{a_{ij}\}$, and with \emph{curving}
described by $\{\omega_i\}$.

\subsection{\v Cech formalism for generalized coverings}
\label{sec:prel_coverings}

In this section, we provide the necessary machinery to translate
statements and computations carried out in a conventional
\v{C}ech covering by open sets to other kinds of coverings, such
as the universal cover, that will allow to merge results from our
previous approach~\cite{prev_paper} into the present one. This
formalism is not yet part of a mathematical physics curriculum,
so here we present the prerequisites necessary for computing
\v{C}ech cohomology, referring to the standard sources
~\cite{artin,milne,sga} where the theoretical background is
explained.

Let $M$ be a smooth manifold or topological space. The general
idea is to pass from inclusions $U\hookrightarrow M$ to general
local homeomorphisms $U\to M$ which are not necessarily
injective. Technically, one fixes a category $\cover{C}_M$ whose
objects are spaces \'etale over $M$, morphisms are the covering
maps, and which is closed with respect to the fiber product of
the maps over $M$, with $M$ being the terminal object in
$\cover{C}_M$.  The coverings are surjective families of local
homeomorphism in $\cover{C}_M$, namely families $\{f_i:U_i\to
U\}$ of $M$-maps such that $U=\bigcup_i f_i(U_i)$. In practice,
we shall restrict our attention to covering maps of $M$ itself.
The key observation is that if $U_i\hookrightarrow M$ and
$U_j\hookrightarrow M$ are inclusions, then $U_i\cap U_j \equiv
U_i\times_M U_j$ --- the fiber product of maps
$U_i\hookrightarrow M$ and $U_j\hookrightarrow M$ --- so that the
notion of fiber product for covering maps replaces the notion of
intersection of open sets.

For a covering $U\to M$ in $\cover{C}_M$ we obtain an augmented
simplicial object~\cite{artin-mazur}
\begin{displaymath}
  \xymatrix{M & U\ar[l]_-\pi & U\times_M U\ar@<-.4ex>[l]
    \ar@<.4ex>[l] & U\times_M U\times_M U\ar@<-.7ex>[l] \ar[l]
    \ar@<.7ex>[l] & \dots \ar@<-1ex>[l]
    \ar@<-.3ex>[l]\ar@<.3ex>[l] \ar@<1ex>[l]}
\end{displaymath}
by considering the nerve $N_\bullet (U\to M)$. Specifically, for
any integer $q\geq 0$ we define
\begin{equation*}
  N_q(U\to M)=\underbrace{U\times_M \dots \times_M
  U}_{(q+1)-\text{times}}
\end{equation*}
where for $i=0,\dots,q$ the arrows are the maps $d_i :N_q(U\to
M)\to N_{q-1}(U\to M)$, forgetting the $i\text{-th}$ factor in
the product.

For an abelian sheaf $\sheaf{F}$ on $M$ (more precisely, on
$\cover{C}_M$) the \v{C}ech complex relative to a covering $U\to
M$ in $\cover{C}_M$ is defined by setting for any $q\geq 0$
\begin{equation*}
  \Check{C}^q(U;\sheaf{F}) = \Gamma (N_q(U\to M),\sheaf{F})
  \quad\text{with}\quad \deltacheck = \sum_{i=0}^q (-1)^i
  d^*_i\,.
\end{equation*}
The ordinary \v{C}ech formalism is recovered by considering an
open cover $\cover{U}_M= \{U_i\}_{i\in I}$ of $M$ and the
covering $\coprod_{i\in I}U_i \to M$, so that in degree $q$ we
just get the disjoint union of all $q$-fold intersections
\begin{equation*}
  N_q(\cover{U}_M) = \coprod_{i_0,\dots ,i_q} U_{i_0}\cap\dots
  \cap U_{i_q}\,,
\end{equation*}
and the resulting \v Cech complex is the standard one.

At the other extreme, let $U\to M$ be a regular covering map and
$G=\mathrm{Deck}(U/M)$ the corresponding group of deck
transformations acting properly on $U$ on the right. One
immediately verifies that
\begin{equation*}
  \underbrace{U\times_M \dots \times_M U}_{(q+1)-\text{times}}
  \cong U\times \underbrace{G\times\dots \times
  G}_{q-\text{times}}\,,
\end{equation*}
and under this isomorphism the maps $d_i :N_q(U\to M) \to
N_{q-1}(U\to M)$ become
\begin{equation*}
  d_i(x,g_1,\dots ,g_q) =
  \begin{cases}
    (x\cdot g_1,g_2,\dots,g_q)& i=0\\
    (x,g_1,\dots,g_ig_{i+1},\dots,g_q) & i=1,\dots,q-1\\
    (x,g_1,\dots,g_{q-1}) & i=q\,.
  \end{cases}
\end{equation*}
Hence, the \v Cech complex with respect to $U\to M$ becomes the
usual Eilenberg-MacLane cochain complex on $G$ with values in the
$G$-module $\sheaf{F}(U)$:
\begin{equation*}
  \Check{C}^q(U;\sheaf{F})\cong C^q(G;\sheaf{F}(U))\,.
\end{equation*}
Thus the \v Cech cohomology of this complex is just the group
cohomology of $G$ with values in the $G$-module $\sheaf{F}(U)$,
where the module structure is given by the pull-back action. A
particular case of special interest for us is when $U$ is the
universal cover of $M$, so that $G=\pi_1(M)$.

The formalism clearly extends to the case where we consider a
complex $\sheaf{A}^\bullet$ of sheaves on $M$ --- typically, the
de Rham complex. The hypercohomology with respect to a covering
$U\to M$ will be the cohomology of the total complex of
$\Check{C}^q(U;\sheaf{A}^p)$.

In some favorable cases, one or both spectral sequences
associated to the double complex above will degenerate at the
first level. Degeneration at the first level of the first
spectral sequence, that is, the one associated to the filtration
on $p$, is equivalent to
\begin{equation*}
  \Check{H}^q(U\to M;\sheaf{A}^p) = 0\quad \text{for all $q>0$.}
\end{equation*}
Since each $\sheaf{A}^p$ is assumed to be a sheaf, that is,
${A}^p(M)$ is the kernel
\begin{equation*}
  \xymatrix{\protect{A}^p(M)\ar[r] &
  \protect\sheaf{A}^p(U)\ar@<.5ex>[r] \ar@<-.5ex>[r] &
  \protect\sheaf{A}^p(U\times_M U)}\,,
\end{equation*}
the cohomology of the total complex $\Check{C}^q(U;\sheaf{A}^p)$
equals $H^p_\mathit{dR}(M)$.

On the other hand, the degeneration of the other spectral
sequence (at the same level) means the complex
$\sheaf{A}^\bullet$ is a resolution of some sheaf $\sheaf{F}$, so
that the total cohomology equals $\Check{H}^p(U\to M;\sheaf{F})$.
Therefore, when both of these cases are realized, we have a
\v{C}ech-de Rham type situation~\cite{bott-tu}, that is
\begin{equation*}
  \HH^p(U\to M;\sheaf{A}^\bullet)\cong \Check{H}^p(U\to
  M;\sheaf{F})\cong H^p_\mathit{dR}(M)\,.
\end{equation*}

The obvious example of this situation is the \v Cech-de Rham
double complex relative to the ordinary cover $\coprod_{i\in
I}U_i$, where the above isomorphism gives the usual de Rham
theorem: $\Check{H}^p(M,\CC)\cong H^p_\mathit{dR}(M)$.  Another
example of utmost importance is the universal cover $\HH\to X$ of
a Riemann surface $X$ of genus $g>1$. Since there exist
$\pi_1(X)$-equivariant partitions of unity~\cite{kra}, the
sheaves $\sheaf{A}^p_\HH$ are acyclic: $H^q(\pi_1(X),
\sheaf{A}^p_\HH)=0$ for $q>0$ and all $p$. Moreover, since $\HH$
is contractible, the de Rham complex $\sheaf{A}^\bullet(\HH)$ is
obviously acyclic in dimension greater than zero, and as a result
we have the isomorphism\footnote{See also~\cite{prev_paper} for a
simple-minded proof without spectral sequences.}
\begin{equation*}
  H^p(\pi_1(X),\CC)\cong H^p_\mathit{dR}(M)\,.
\end{equation*}

\subsection{Evaluation over the fundamental class}
\label{sec:prel_fundamental}
For the construction of the action functional we need to evaluate
Deligne cohomology classes against the fundamental class $[X]$ of
$X$, which we need to represent as a cycle in a suitable
homological double complex --- in a way analogous to the use of
\v Cech resolutions to compute the hypercohomology.

The aim of this section is to introduce the minimum set of tools
necessary to describe the homological (double) complex and to
perform the evaluation, relegating all technical details to the
appendix. There, we construct an explicit representative $\Sigma$
of $[X]$ with respect to a covering $U\to X$ by mirroring the
cohomology computations done in~\ref{sec:local_cocycle}.  The
computations are explicit enough that the reader who is only
interested in the formulas for $\Sigma$ can
read~\ref{sec:fund_class_computation} directly.  Also, the reader
interested only in the construction of the local action cocycle
can safely proceed to sect.~\ref{sec:Action}.

As usual, whenever we mention facts that are not specific to $X$
being a Riemann or topological surface, we use the notation $M$
to denote a general smooth manifold or topological space with
covering $U\to M$.

\subsubsection{}
Consider the double complex
\begin{equation*}
  \SSS_{p,q}=S_p(N_q(U\to M))\,,
\end{equation*}
where $N_\bullet(U\to M)$ is the nerve of the covering $U\to M$
and $S_\bullet$ is the singular simplices functor, i.e., $S_p(M)$
is the set of continuous maps $\Delta^p\to M$, where $\Delta^p$
is the standard simplex. For every $p\geq 0$, the covering map
$U\to M$ induces a corresponding map $\epsilon :
\SSS_{p,0}=S_p(U)\to S_p(M)$ between simplices --- the
augmentation map. The double complex $\SSS_{\bullet,\bullet}$ has
two boundary operators: the usual boundary operator on singular
chains, $\delp:\SSS_{p,q}\rightarrow \SSS_{p-1,q}$, and the
boundary operator $\delpp :\SSS_{p,q}\rightarrow \SSS_{p,q-1}$
induced by the face maps of the nerve: $\delpp=\sum (-1)^i
{d_i}_\star$, where $d_i:N_q(U)\to N_{q-1}(U)$ and ${d_i}_\star$
is the induced map on singular chains. As usual, we have the
simple complex $\Tot\SSS$ with total differential $\del = \delp +
(-1)^p\delpp$ on $\SSS_{p,q}$.

If $U$ is the ordinary \v{C}ech covering
$\cover{U}_M=\coprod_{i\in I}U_i$, then
\begin{equation*}
  S_p(N_q(\cover{U}_M)) = \prod_{i_0,\dots ,i_q}
  S_p(U_{i_0}\cap\dots \cap U_{i_q})\,.
\end{equation*}
If, on the other hand, $U$ is a regular covering space with $G$
as group of deck transformations, then $S_p(U)$ is a $G$-module
with $G$-action given by translation of simplices. It follows
that $S_p(N_q(U))$, for $q>0$, consists of simplices into $U$
parameterized by $q$-tuples of elements in $G$. Taking into
account the expression for the face maps $d_i$, computed
in~\ref{sec:prel_coverings}, we get
\begin{equation*}
  S_p(N_q(U))= S_p(U)\otimes_{\ZZ G}B_q(G)\,,
\end{equation*}
where $B_\bullet(G)$ is the bar resolution~\cite{maclane} and
$\ZZ G$ is the integral group ring of $G$. Hence, for any $p$,
the $\delpp$-homology is just the group homology
\begin{equation*}
  H_q(S_p(N_\bullet(U)))=H_q(G;S_p(U))\,.
\end{equation*}

We are interested in the situation when $\SSS_{\bullet,\bullet}$
has no homology with respect to the second index, except in
degree zero, namely we want
\begin{equation*}
  H_q(S_p(N_\bullet(U\to M)))\cong
  \begin{cases}
    S_p(M)& q=0\\ 0 & q >0\,,
  \end{cases}
\end{equation*}
for the $\delpp$ homology. In this case we say that
$\SSS_{p,\bullet}$ \emph{resolves} $S_p(M)$ and one has the
isomorphism
\begin{equation*}
  H_\bullet(M,\ZZ)\equiv H_\bullet(S_\bullet(M)) \cong
  H_\bullet(\Tot \SSS_{\bullet,\bullet})\,.
\end{equation*}
This isomorphism is induced by the augmentation map
$\epsilon:\Tot\SSS \to S_\bullet(M)$, which assigns to any chain
$\Sigma$ of total degree $n$ in $\Tot\SSS$ the chain
$\epsilon(\Sigma_{n,0})$, where $\Sigma_{n,0}$ is the component
in $\SSS_{n,0}$. It is easy to see that this map is a chain map,
it sends cycles into cycles and induces the above isomorphism.
Details can be found, e.g., in~\cite{maclane}.  \footnote{A
  detailed calculation along these lines can also be found in the
  appendix of~\protect\cite{prev_paper}.} Observe that this
situation is realized for both the examples of an open \v{C}ech
cover and of a regular covering $U\to M$ (cf.~the appendix). For
completeness, in the appendix we briefly analyze the implications
of the requirement that the double complex
$\SSS_{\bullet,\bullet}$ is acyclic with respect to the first
index, and their relations with good covers.

\subsubsection{}
\label{sec:sigma}
For a topological manifold $M$ of dimension $n$, we need to
represent $[M]$ with a total cycle $\Sigma$ of degree $n$ in
$\Tot\SSS_{\bullet,\bullet}$. It has the form
\begin{displaymath}
\Sigma = \Sigma_0 + \sum_{k=1}^n\,
(-1)^{\sum_{l=0}^{k-1}(n-l)}\Sigma_k\,,
\end{displaymath}
where $\Sigma_k\in \SSS_{n-k,k}$ and
\begin{equation*}
\delp \Sigma_0 = \delpp \Sigma_1\,,\dots\,, \delp \Sigma_{k-1} =
\delpp \Sigma_k\,,\dots\,, \delp \Sigma_n=0\,.
\end{equation*}
The choice of signs ensures $\del \Sigma=0$, where $\del$ is the
total differential in $\Tot\SSS_{\bullet,\bullet}$.  By
definition, $\Sigma$ is a ``lift'' of $M$ considered as a chain
in $S_n(M)$, i.e.~$\epsilon(\Sigma_0)= M$, where
$M=\sum_i\sigma_i$ for a suitable collection of singular
simplices $\sigma_i\in S_n(M)$. The existence of the elements
$\Sigma_1\,, \dots\,, \Sigma_n$ follows from the
$\delpp$-exactness assumption and the fact that $\Sigma_0$ lifts
$M$. Indeed, we have $0=\del M =\del\epsilon (\Sigma_0) =
\epsilon(\delp \Sigma_0)$, so that there exists $\Sigma_1\in
\SSS_{n-1,1}$ such that $\delp\Sigma_0=\delpp\Sigma_1$, and so
on.

Specializing to the case when $M\equiv X$ is a Riemann surface,
the representative of the fundamental class $[X]$ is the cycle
$\Sigma=\Sigma_0 +\Sigma_1 -\Sigma_2$, with components
$\Sigma_k\in\SSS_{2-k,k}$ satisfying $\delp\Sigma_0 =
\delpp\Sigma_1$, $\delp\Sigma_1= \delpp\Sigma_2$, and
$\delpp\Sigma_0=\delp\Sigma_2=0$. This cycle is explicitly
constructed in the appendix for the case of an ordinary \v{C}ech
cover $U=\cover{U}_X$ and in~\cite{prev_paper} for the case of
the universal cover $\HH\to X$.

Here we present the basic formulas for the \v{C}ech case, which
also gives the flavor of the general procedure which carries over
to the other coverings unchanged.

Following~\cite{goldberg,weil}, introduce the symbol
$\Delta_{i_0,\dots,i_q}$ to denote the $(q+1)$-fold intersection
thought as a generator in $\SSS_{p,q}$, so that a generic element
can be written in the form:
\begin{displaymath}
\sigma = \sum_{i_0\dots i_q}\sigma_{i_0\dots i_q}\cdot
\Delta_{i_0\dots i_q}\,,
\end{displaymath}
where $\sigma_{i_0\dots i_q}$ are are $p$-simplices for
$U_{i_0}\cap\dots \cap U_{i_q}$, i.e.~continuous maps
$\Delta^p\to U_{i_0}\cap\dots \cap U_{i_q}$.  It is immediate to
verify that
\begin{displaymath}
\delpp\Delta_{i_0\dots i_q} = \sum_{j=0}^q (-1)^j
\Delta_{i_0,\dots ,\widehat{i_j},\dots ,i_q}\,,
\end{displaymath}
where the $\,\widehat{}\,$ sign denotes omission. Then
\begin{displaymath}
\delpp \sigma = \sum_{i_0,\dots ,i_{q-1}} \bigl (
\sum_{k=0}^q\sum_j (-1)^k\sigma_{i_0,\dots,
\underset{\overset{\uparrow}{k\text{-th}}}{j}, \dots ,i_q}
\bigr)\cdot \Delta_{i_0,\dots ,i_{q-1}},
\end{displaymath}
where the summation goes over ordered sets of indices (it is
assumed that $I$ is an ordered set). Thus with the convention
that $\sum_{\langle i_0,\dots,i_q\rangle}$ is the sum over sets
of indices $\{i_0,\dots ,i_q\}$ with $i_0\leq\dots\leq i_q$, we
can rewrite the last equation as
\begin{equation*}
\delpp\sigma=\sum_{k=0}^q (-1)^k \sum_{\langle
i_0,\dots,\underset{\overset{\uparrow}{k\text{-th}}}{j},\dots
,i_q\rangle} \sigma_{i_0,\dots,j,\dots ,i_q}\cdot
\Delta_{i_0,\dots ,i_{q-1}} \,.
\end{equation*}

Now, consider the problem of constructing the total cycle
$\Sigma=\Sigma_0 + \Sigma_1-\Sigma_2$ representing $[X]$.
Representing the components $\Sigma_i$ as:
\begin{displaymath}
\Sigma_0=\sum_i\sigma_i\cdot \Delta_i\,, \qquad
\Sigma_1=\sum_{\langle ij\rangle}\sigma_{ij}\cdot \Delta_{ij}\,,
\qquad \Sigma_2=\sum_{\langle ijk\rangle}\sigma_{ijk}\cdot
\Delta_{ijk}\,,
\end{displaymath}
we first construct $\Sigma_0$ as follows. Starting from the
nerve of the cover $\cover{U}_X$ consider a triangulation of $X$
by $\cover{U}$-small simplices, i.e.~each simplex comprising the
triangulation has support in some open set $U_i$ belonging to the
cover (cf.~the appendix for the detailed procedure). Then
$X=\sum_i \sigma_i$, where each chain $\sigma_i$ is a sum of
simplices whose support is contained in $U_i$ for each $i$, and
one immediately writes $\Sigma_0=\sum_i\sigma_i\cdot \Delta_i$.
The other components are determined by the $\delpp$-exactness
condition of the complex. Namely, from the above expression and
$\delp\Sigma_0 = \delpp\Sigma_1$, $\delp\Sigma_1= \delpp\Sigma_2$
one gets the equations
\begin{gather*}
\sum_i\delp\sigma_i\cdot \Delta_i = \bigl(\sum_{\langle
i,j\rangle}\sigma_{ji} -\sum_{\langle i,j\rangle}\sigma_{ij}
\bigr)\cdot \Delta_i\,,\\ \sum_{\langle i,j
\rangle}\delp\sigma_{ij}\cdot \Delta_{ij} = \bigl(\sum_{\langle
k,i,j\rangle}\sigma_{kij} -\sum_{\langle
i,k,j\rangle}\sigma_{ikj} +\sum_{\langle
i,j,k\rangle}\sigma_{ijk}\bigr)\cdot \Delta_{ijk}
\end{gather*}
for components $\sigma_{ij}$ and $\sigma_{ijk}$. Explicit
expression for these components in terms of the barycentric
decomposition is given in the appendix.

\subsubsection{}
\label{sec:pairing}

In order to discuss the evaluation pairing, we need to address
the issue of the index shift in the Deligne complex. One way is
to explicitly use the exponential map described
in~\ref{sec:prel_sheaves} to revert the indexing to the familiar
form without a shift, at the cost of introducing an explicit
multiplicative structure via the exponential. Another way is to
introduce an \emph{ad hoc} index shift in homology to mirror the
one in the Deligne complex, i.e.~to consider singular
$q$-simplices to be of homological degree $q+1$. The resulting
pairing will be additive, but only defined $\mod\ZZ(p)$. The two
approaches are in the end the same.

We start with second approach. Let $(\KKK_\bullet,\del)$ be a
\emph{homological} complex. The canonical way to shift it is to
introduce $\KKK[1]_\bullet \eqdef\KKK_{\bullet -1}$, with
$\del_{[1]}=-\del$, cf.~\cite{maclane}. We require instead that
the new differential be simply $\del$, while retaining the index
shift. Thus we replace $\SSS_{r,s}=S_r(N_s(U\to M))$ by the new
double complex
\begin{displaymath}
  {}^\backprime\SSS_{r,s}=S_{r-1}(N_s(U\to M))\,,
\end{displaymath}
with differential $\del = \delp + (-1)^r\delpp$, where $\delp$ is
the usual singular boundary, as before. If $\Sigma =
(\Sigma_0,\cdots,\Sigma_q)$, with $\Sigma_k\in \SSS_{q-k,k}$, is
a $q$-chain in $\Tot\SSS_{\bullet,\bullet}$, then it maps to the
$(q+1)$-chain ${}^\backprime\Sigma = ((-1)^q\Sigma_0,
\cdots,(-1)^{q-k}\Sigma_k,\cdots,\Sigma_q,0)$ in
$\Tot{}^\backprime\SSS_{\bullet,\bullet}$, and
${}^\backprime\Sigma$ is a cycle if and only if $\Sigma$ is a
cycle.

Let $\CCC^{\bullet,\bullet}$ be a \v{C}ech resolution of the
Deligne complex $\ZZ(p)^\bullet_\mathcal{D}$ with respect to the
covering $U\to M$. The pairing between $\CCC^{r,s}$ and
${}^\backprime\SSS_{r,s}$ is defined as follows
(cf.~\cite{goldberg,weil}). To the pair $(\phi,\sigma)$, where
$\phi$ is a collection $\{\phi_{i_0,\dots,i_s}\}$ of
$(r-1)$-forms on $N_s(U\to M)$ for $r>0$, or integers $\ZZ(p)$
for $r=0$ and $\sigma=\sum\sigma_{i_0,\dots,i_s}\cdot
\Delta_{i_0,\dots,i_s}\in {}^\backprime\SSS_{r,s}$ we assign
\begin{equation}\label{pairing}
  \langle\phi,\sigma\rangle =
  \begin{cases}
  \sum_{\langle i_0,\dots,i_s\rangle}
  \int_{\sigma_{i_0,\dots,i_s}}\phi_{i_0,\dots,i_s} & r > 0 \\ 0
  & r =0\,.
  \end{cases}
\end{equation}
To extend this pairing to $\Tot \CCC^{\bullet,\bullet}$ and
$\Tot{}^\backprime\SSS_{\bullet,\bullet}$, let
${}^\backprime\sigma =(\sigma_0,\sigma_1,\cdots,\sigma_{n-1},0)$,
with $\sigma_k\in {}^\backprime\SSS_{n-k,k}$, and $\Phi =
(\phi_0,\phi_1,\cdots,\phi_n)$, with $\phi_k\in
\CCC^{n-k,k}$.  Then we define
\begin{equation}
  \label{eq:pairing_total}
  \langle\Phi,{}^\backprime\sigma\rangle =
  \sum_{k=0}^{n-1}\langle\phi_k,\sigma_k\rangle\,,
\end{equation}
where, $\phi_k=0$ for all $k<n-p$, if, of course, $n>p$. Note
that so far the pairing was defined to have values in $\CC$.
However, the fundamental fact is that \emph{away from the
  truncation degree}, i.e.~when the total degree $n$ is strictly
less than $p$, and therefore the form degree is strictly less
than $p-1$, the total differentials $D$ and $\del$ are transpose
to each other modulo $\ZZ(p)$:
\begin{equation}
  \label{dual}
  \langle D \Phi,{}^\backprime\Sigma\,\rangle =\langle\,\Phi,\del
  {}^\backprime\Sigma \rangle \mod \ZZ(p)\,.
\end{equation}
for all $\Phi \in \Tot\CCC^{\bullet,\bullet}$,
${}^\backprime\Sigma \in
\Tot{}^\backprime\SSS_{\bullet,\bullet}$.  This readily follows
from the very definition of the Deligne complex. Equation
\eqref{dual} means that the pairing $\langle\;,\;\rangle$,
\emph{considered modulo} $\ZZ(p)$, defines a pairing between
$H^{\bullet}(\Tot\CCC^{\bullet,\bullet})$ and
$H_{\bullet}(\Tot\SSS_{\bullet,\bullet})$ away from the
truncation degree $p-1$.

Formula \eqref{dual} would not hold for degrees bigger or equal
than $p-1$, \emph{unless} $\dim M=p-1$ --- the case where the
truncation becomes unimportant. This is the situation we will be
interested in in sect.~\ref{sec:Action}. Therefore in this case
the pairing \eqref{eq:pairing_total} descends to the
corresponding homology and cohomology groups and is non
degenerate. It defines a pairing between
$H^{\bullet}(\Tot\CCC^{\bullet,\bullet})$ and
$H_{\bullet}(\Tot\SSS_{\bullet,\bullet})$ which we continue to
denote by $\langle\;,\;\rangle$.

Let us show how these formulas work in the case of a Riemann
surface $X$ and a Deligne cocycle $\Omega =
(\omega_i,a_{ij},f_{ijk},m_{ijkl})$ of total degree $3$. (Recall
that the individual elements are subjects to the
relations~\eqref{eq:del_3_cocycle}.) Let $\Sigma =
(\Sigma_0,\Sigma_1,-\Sigma_2)$ be a representative in
$\Tot\SSS_{\bullet,\bullet}$ of the fundamental class $[X]$. Then
the corresponding element in the shifted complex will be
\begin{displaymath}
{}^\backprime\Sigma = (\Sigma_0,-\Sigma_1,-\Sigma_2,0)\,.
\end{displaymath}
Omitting the indices, the evaluation of the class of $\Omega$
over $[X]$ will be computed by the expression
\begin{equation}
\label{eq:evaluation}
\langle\Omega,{}^\backprime\Sigma\rangle = \langle
\omega,\Sigma_0\rangle -\langle a, \Sigma_1\rangle -\langle
f,\Sigma_2 \rangle\,,
\end{equation}
where each term in~\eqref{eq:evaluation} should be expanded
according to~\eqref{pairing}. This evaluation takes its values in
$\CC/\ZZ(3)$ and does not depend on the representative cocycle of
the Deligne cohomology class $[\Omega]\in
H^3_\mathcal{D}(M,\ZZ(3))$.

Another way to define the pairing is to use explicitly the
quasi-isomorphism
\begin{displaymath}
  \ZZ(p)^\bullet_\mathcal{D}\cong (\sheaf{A}^*_M
  \xrightarrow{d\log}\sheaf{A}^1_M \overset{d} {\longrightarrow}
  \dots \overset{d}{\longrightarrow} \sheaf{A}^{p-1}_M)[-1]\,,
\end{displaymath}
induced by the exponential map (see ~\ref{sec:prel_sheaves}).  In
this way a cocycle representing a class of degree $k$ in
$H^k_\mathcal{D}(M,\ZZ(p))$ becomes a cocycle of degree $k-1$ in the double
complex
\begin{equation*}
 \Check{C}^\bullet(U\to M, \sheaf{A}^*_M \to \sheaf{A}^1_M \to
 \dots\to\sheaf{A}^{p-1}_M) \,.
\end{equation*}
In particular, if $\Omega=(\omega_i,a_{ij},f_{ijk},m_{ijkl})$, subject 
to the
relations \eqref{eq:del_3_cocycle}, is a cocycle of total degree
3 in $\Check{C}^\bullet(U\to M, \ZZ(3)^\bullet_\mathcal{D})$
representing a Deligne class of total degree 3, the element
\begin{equation*}
  \Psi=\Biggl(\frac{1}{(2\pi i)^2}\omega_i,-\frac{1}{(2\pi
    i)^2}a_{ij},\exp\biggl(\frac{f_{ijk}}{(2\pi i)^2}\biggr)
    \Biggr)
\end{equation*}
is the corresponding cocycle of total degree 2.

As in the previous discussion, we will consider only the the case
when $\dim M=p-1$, where $p$ is the length of the Deligne
complex, so that the truncation becomes irrelevant.  Denote by
$\Tilde{\CCC}^{\bullet,\bullet}$ the double complex
$\Check{C}^\bullet(U\to M,\sheaf{A}^*_M \to \sheaf{A}^1_M \to
\dots\to\sheaf{A}^{n}_M)$, with $n=p-1$. 

Then there exists a natural pairing between $\Tilde{\CCC}^{r,s}$
and $\SSS_{r,s}$ which assigns to the pair $(\psi,c)$ the
evaluation of the $r$-form $\psi$ over a chain
$c=\sum\sigma_{i_0,\dots,i_s}\cdot\Delta_{i_0,\dots,i_s} \in
S_r(N_s(U\to M))$:
\begin{equation*}
\langle\psi,c\rangle =
\int_{\sigma_{i_0,\dots,i_s}}\psi_{i_0,\dots,i_s}\,,
\end{equation*}
with the understanding that for $r=0$ this is just the pointwise
evaluation of an invertible function, defined through the
exponential map. To define a multiplicative pairing between $\Tot
\Tilde{\CCC}^{\bullet,\bullet}$ and $\Tot\SSS_{\bullet,\bullet}$,
let $C=(c_0,c_1,\dots,c_n)$, with $c_i\in \SSS_{n-i,i}$, and
$\Psi = (\psi_0,\psi_1,\dots,\psi_n)$, with $\psi_i\in
\Tilde{\CCC}^{n-i,i}$. Then we define
\begin{equation}\label{eq:pairing_total_mult}
  \langle\Psi,C\rangle_m = \prod_{i=0}^{n-1}
  \exp(\langle\psi_i,c_i\rangle)\cdot \langle
  \psi_n,c_n\rangle\,\in\CC^*.
\end{equation}
By the very construction of the double complexes
$\Tilde{\CCC}^{\bullet,\bullet}$ and $\SSS_{\bullet,\bullet}$,
the total differentials $D$ and $\del$ are transpose to each
other, namely
\begin{equation*}
  \langle D \Psi,C\,\rangle_m =\langle\,\Psi,\del C \rangle_m
\end{equation*}
for all $\Psi \in \Tilde{\CCC}^{\bullet,\bullet}$, $C \in
\SSS_{\bullet,\bullet}$. The pairing
\eqref{eq:pairing_total_mult} descends to the corresponding
homology and cohomology groups and is non degenerate. It defines
a pairing between
$H^{\bullet}(\Tot\Tilde{\CCC}^{\bullet,\bullet})$ and
$H_{\bullet}(\Tot\SSS_{\bullet,\bullet})$ which we continue to
denote by $\langle\;,\;\rangle_m$.

It is easy to describe the relation between the multiplicative
pairing $\langle\;,\;\rangle_m$ and the $\CC/\ZZ(p)$-valued
additive pairing introduced earlier. Namely, let
$\Phi\in\Tot\CCC^{\bullet,\bullet})$,
$C\in\SSS_{\bullet,\bullet}$ and let $\Psi$ be the corresponding
element in $\Tilde{\CCC}^{\bullet,\bullet}$. The we have
\begin{equation*}
\langle\Psi,C\,\rangle_m =\exp\{\langle\Phi,^\backprime C\,
\rangle/(2\pi i)^{p-1}\}.
\end{equation*}

It what follows we will use freely both forms of the pairing,
multiplicative and additive, depending on the context.

\section{Construction of the action}
\label{sec:Action}
\subsection{General remarks}
The next sections will be devoted to the detailed construction of
the action functional --- or rather its exponential --- by
specifying the following.
\begin{itemize}
\item[a.] A resolution of the Deligne complex
  $\ZZ(3)^\bullet_\mathcal{D}$.
\item[b.] A representative for a class in
  $H^3_\mathcal{D}(X,\ZZ(3))$ that ``starts'' from a collection
  $\{\omega_i[f]\}_{i\in I}$ of "local Lagrangians densities''
  --- top forms on $X$ --- defined with respect to a given
  covering $\cover{U}_X =\{U_i\}_{i\in I}$ of $X$.
\end{itemize}
The latter data come from Polyakov's ansatz, with dynamical field
given by a deformation map $f:X\rightarrow \Tilde X$ and with
external field given by a smooth projective connection of $X$.
Before doing so, we make some remarks of general character.
\begin{itemize}
\item The Deligne complex $\ZZ(3)^\bullet_\mathcal{D}$ is
  especially convenient for treating various logarithmic terms
  produced in descent calculations, while keeping additivity.
\item The ``local Lagrangian'' $\Omega[f]$ appears as a total
  cocycle of total degree 3 in the Deligne complex
  $\ZZ(3)^\bullet_\mathcal{D}$, and we define the action
  functional by evaluating this cocycle over the representative
  $\Sigma$ of the fundamental class of the Riemann surface $X$
  \begin{displaymath}
  S[f]=\langle\,\Omega[f],\Sigma \rangle,
  \end{displaymath}
  described in~\ref{sec:sigma}. According to~\ref{sec:pairing},
  $S[f]\in\CC/\ZZ(3)$, so that the functional
  \begin{displaymath}
  A[f]=\langle\,\Omega[f],\Sigma \rangle_m=\exp\{S[f]/(2 \pi
  i)^2\}
  \end{displaymath}
  is the \emph{exponential} of the action.
\item A similar approach was taken in~\cite{alvarez,gawedzki} in
  order to describe certain topological terms arising in
  two-dimensional quantum field theories. In our case the field
  is a deformation $f:X\rightarrow \Tilde X$ and the procedure
  differs in that we construct the whole representing cocycle
  starting from one end of the descent staircase.
\item According to~\cite{freed,del-freed} the exponentials of
  action functionals should be more properly regarded as
  $\CC^*$-torsors rather than numbers. This is most apparent when
  dealing with manifolds with boundaries. A similar situation
  occurs in our case, when $X$ is a compact Riemann surface: the
  definition of the local Lagrangian cocycle $\Omega[f]$ depends
  on the trivialization of the tame symbol $(TX,TX]$, described
  by an ($f$-independent) element of $H^2(X,\CC^*)\cong \CC^*$.
  As a result, the multiplicative action functional $A[f]$ is a
  $\CC^*$-torsor.
\item The action functional $A[f]$, defined through
  hypercohomology admits the following geometric interpretation.
  According to Sect.~\ref{sec:prel_sheaves}, the group
  $H^3_\mathcal{D}(X,\ZZ(3))$ classifies isomorphism classes of
  gerbes equipped with connective structure and
  curving~\cite{bry_mcl1,bry_book}. Since $\dim X =2$, these are
  necessarily flat, therefore they are classified by their
  holonomy via the isomorphism $H^3_\mathcal{D}(X,\ZZ(3))\cong
  H^2(X, \CC^*)$. Thus $A[f]$ can be interpreted as the holonomy
  of an appropriate higher algebraic structure.
\end{itemize}

\subsection{Setup for regular \v Cech coverings}
\label{sec:action_cech}

Let $\cover{U}_X=\{U_i\}_{i\in I}$ be an open cover of $X$,
which we assume to be a good cover, i.e.~all nonempty
intersections $U_{i_0,\dots ,i_p}=U_{i_0}\cap \dots \cap U_{i_p}$
are contractible. Therefore, we are in a \v Cech-de Rham
situation~\cite{bott-tu,weil}, and the double complex
$\CCC^{p,q}_\mathcal{D}\eqdef
\Check{C}^q(\cover{U}_X,\ZZ(3)^p_\mathcal{D})$ computes
$H^\bullet_\mathcal{D}(X,\ZZ(3))$. Let $\{z_i:U_i \rightarrow
\CC\}_{i\in I}$ be holomorphic coordinates for the complex
structure of $X$, and let $z_{ij}: z_j(U_i\cap U_j)\rightarrow
z_i(U_i\cap U_j)$ be coordinate change functions:
$z_i=z_{ij}\circ z_j$ on $U_i\cap U_j$.
\begin{remark}
  One could also use coordinate functions with in $\PP^1$ instead
  of $\CC$.
\end{remark}

More generally, for $U_{i_0,\dots ,i_q}$ there are holomorphic
coordinates $z_{i_0},\dots ,z_{i_q}$ with $z_{i_k}=z_{i_k
  i_{k+1}}(z_{i_{k+1}})$, $k=0,\dots,q-1$. If
$\phi\in\Check{C}^{q-1}(\cover{U}_X,\sheaf{A}^p_X)$ is a \v Cech
cochain, i.e.~ $\phi=\{\phi_{i_0,\dots,i_{q-1}}\}$, where the
components $\phi_{i_0,\dots,i_{q-1}}$ are $p$ forms on
$U_{i_0,\dots ,i_{q-1}}$ its \v Cech differential is defined as
\begin{equation*}
  \deltacheck\phi_{i_0,\dots,i_q} = \sum_{k=0}^{q-1} (-1)^k
  \phi_{i_0,\dots,\Hat{i_k},\dots,i_q} +(-1)^q
  (z_{i_{q-1}i_q})^*\phi_{i_0,\dots,i_{q-1}}\,.
\end{equation*}
It is understood that each component $\phi_{i_0,\dots, i_{q-1}}$
of a \v Cech cochain $\phi$ 
is expressed in the coordinate $z_{i_{q-1}}$, i.e.~the one
determined by the last index, and we will use this convention
throughout the paper.

Given a quasi-conformal map $f:X\to \Tilde X$, denote by
$\cover{V}_X=\{V_i\}_{i\in I}$, where $V_i=f(U_i)$, the
corresponding good open cover for $\Tilde X$. Let $\{w_i:V_i
\rightarrow \CC\}_{i\in I}$ be holomorphic coordinates for the
complex structure of $\Tilde X$, and let $w_{ij}: w_j(V_i\cap
V_j) \rightarrow w_i(V_i\cap V_j)$ be the corresponding
coordinate change functions: $w_i=w_{ij}\circ w_j$ on $V_i\cap
V_j$. Let $f_i=w_i\circ f\vert_{U_i}\circ z^{-1}_i,\,i\in I$, be
local representatives of the map $f$, satisfying the
transformation law
\begin{equation}\label{eq:f_transf}
  f_i\circ z_{ij} = w_{ij}\circ f_j\,.
\end{equation}

Denote $\del f_i\eqdef \del f_i/\del z_i$ and $\delb_i f_i\equiv
\delb f_i\eqdef \del f_i/\del \bar z_i$, and introduce local
representatives of the Beltrami differential $\mu$ by $\mu_i=
\delb f_i/\del f_i$.

It follows from \eqref{eq:f_transf} that
\begin{gather}
  \label{eq:delf_transf}
  \del f_i\circ z_{ij}\cdot z^\prime_{ij} = w^\prime_{ij}\circ
  f_j\cdot \del f_j, \\
  \label{eq:delbf_transf}
  \delb f_i\circ z_{ij}\cdot \overline{z^\prime_{ij}} =
  w^\prime_{ij}\circ f_j\cdot\delb f_j, \\ \intertext{and}
 \label{eq:m_transf}
 \mu_i\circ z_{ij}\cdot
  \frac{\overline{z^\prime_{ij}}}{z^\prime_{ij}} = \mu_j\,.
\end{gather}

Since $\xi_{ij}\eqdef z^\prime_{ij}\circ z_j=dz_i/dz_j$ are
transition functions for the holomorphic tangent bundle $TX$, and
$\tilde\xi_{ij}\eqdef w^\prime_{ij}\circ w_j= dw_i/dw_j$ are the
transition functions for $T\tilde X$, it follows from
\eqref{eq:delf_transf} that $\del f$ is a section of the bundle
$T^*X\otimes f^{-1}T\Tilde X$, or $\del f \in
A^{1,0}(f^{-1}T\Tilde X)$. Here $f^{-1}T\tilde X$ is the
pull-back of the tangent bundle over $\tilde X$ by
$f$. Similarly, $\delb f \in A^{0,1}(f^{-1}T\Tilde X)$.

In the $C^\infty$ category $f^{-1}T\Tilde X\cong TX$, so that
$T^*X\otimes f^{-1}T\Tilde X$ is isomorphic to the trivial
bundle.  This is also implied directly by the transition formula
\eqref{eq:delf_transf}, since $\del f_i\neq 0$, $f$ being a
diffeomorphism. Thus $\del f$ is an explicit trivializing section
for $T^*X\otimes f^{-1}T\Tilde X$, that establishes the
isomorphism between $T^*X\otimes f^{-1}T\Tilde X$ and the trivial
line bundle.

Introducing representatives $c_{ijk}$ and $\Tilde c_{ijk}$ for
the first Chern classes $c_1(TX)=c_1(\Tilde TX)$, we have
\begin{subequations}
  \label{eq:branches}
  \begin{align}
    c_{ijk} &= \deltacheck (\{\log
    z^\prime_{\cdot\cdot}\})_{ijk}\,,\\ \Tilde c_{ijk} &=
    \deltacheck (\{\log w^\prime_{\cdot\cdot}\circ
    f_\cdot\})_{ijk}\,,\\ b_{ij} &= \log w^\prime_{ij}\circ f_j -
    \log z^\prime_{ij} - \log \del f_i\circ z_{ij} + \log \del
    f_j \,,
  \end{align}
\end{subequations}
and, obviously, $\deltacheck (\{b_{\cdot\cdot}\})_{ijk} =
\Tilde{c}_{ijk} - c_{ijk}$. All the numbers $b_{ij}$, $c_{ijk}$
and $\Tilde c_{ijk}$ are in $\ZZ(1)$.

Although one can get $\Tilde c_{ijk}=c_{ijk}$ and $b_{ij}=0$
through a suitable redefinition of the logarithm branches, there
is no additional complication (except, perhaps, the notation) in
keeping the general situation.

\subsection{The local Lagrangian cocycle}
\label{sec:local_cocycle}

In order to construct the action functional, one needs an ansatz
for its top degree part. Following~\cite{prev_paper}, we promote
the standard Polyakov's chiral action\footnote{More precisely,
  Polyakov's chiral action has no second term in
  \eqref{eq:omega}, which, in fact, is not necessary in genus
  zero.},
\begin{equation}
  \label{eq:omega}
  \omega_i=\frac{\del^2f_i}{\del f_i}\, \del\mu_i\, dz_i\wedge
  d\Bar z_i +2\,\mu_i\,h_i\,dz_i \wedge d\Bar z_i \,,
\end{equation}
to an element $\{2\pi\sqrt{-1}\omega_i\}_{i\in I}\in
\CCC^{3,0}_\mathcal{D}$. Here $h=\{h_i\}_{i\in I}$ is a
$C^\infty$ coboundary for the Schwarzian cocycle
\begin{displaymath}
  \{z_i,z_j\}
  =\frac{d^3z_i}{dz^3_j}
  -\frac{3}{2}\Bigl(\frac{d^2z_i}{dz^2_j}\Bigr)^2,
\end{displaymath}
relative to the cover $\cover{U}_X$ (see~\cite{gunning}). In
other words, it satisfies the following transformation law
\begin{equation}
  \label{eq:h}
  \{z_{i},z_j\}=h_j - h_i\circ z_{ij}\cdot (z^\prime_{ij})^2\,
\end{equation}
on $U_i\cap U_j$.  Clearly, such an $h$ exists, since the
Schwarzian cocycle is already zero in the holomorphic
category~\cite{gunning}. The space $\mathcal{Q}(X)$ of all such
$h$ includes the holomorphic projective connections, and is an
affine space over the vector space
$H^0(X,(\sheaf{A}^{1,0}_X)^{\otimes 2})$. Let us call such an $h$
a \emph{smooth} projective connection (even though that we do not
relate it to projective structures).

Following the usual strategy~\cite{FS} of descending the
staircase in the double complex
$\CCC^{\bullet,\bullet}_\mathcal{D}$, starting with the 0-cochain
$\{\omega_i\}$ of 2-forms on $X$, we find a 1-cochain of 1-forms
$\{\theta_{ij}\}$ and a 2-cochain of functions $\{\Theta_{ijk}\}$
satisfying
\begin{align*}
  \deltacheck(\omega_\cdot)_{ij} &= d\theta_{ij}\\
  \deltacheck(\theta_{\cdot\cdot})_{ijk} &= d\Theta_{ijk}.
\end{align*}
Imposing the condition $\deltacheck{\Theta}=0\mod \ZZ(2)$ ensures
that the total element
\begin{equation*}
\Omega\eqdef 2\pi\sqrt{-1}\bigl(\{\omega_i\}, \{\theta_{ij}\},
  \{-\Theta_{ijk}\}, \{-m_{ijkl}\}\bigr)
\end{equation*}
where $m=\deltacheck\Theta$, is a cocycle in the total complex.

Solvability of the descent equations is proved in the standard
way using the acyclic property of the good cover $\cover{U}_X$
and Poincar\'e lemma on differential forms. Namely, $\deltacheck
d \omega =0$ implies $\deltacheck\omega = d\theta$ and
$0=\deltacheck d\theta = d\deltacheck \theta$ implies
$\deltacheck\theta = d\Theta$.  Finally, from $\deltacheck
d\Theta= d\deltacheck\Theta=0$ one concludes
$\deltacheck\Theta\in \Check{Z}^3(\cover{U}_X,\CC_X)$. {}From de
Rham theorem $\Check{H}^p(X,\CC)\cong H^{p}_{dR}(X)$ it follows
for dimensional reasons that $\deltacheck\Theta = 0$, after
possible rescaling of constants.

The foregoing shows that one can get a ``minimal'' cocycle with
the condition $m_{ijkl}=0$, albeit not in explicit form.
However, our goal is to have a cocycle $\Omega[f]$ with ''good''
dependence on the dynamical field $f$ (i.e.~with the same
variational properties as in the genus zero case).  It is most
remarkable that such cocycle $\Omega[f]$ can in fact be computed
explicitly, allowing for a geometric interpretation as to why
$\deltacheck\Theta =0 \mod \ZZ(2)$. This computation is
accomplished in the following steps.

\subsubsection{$\protect\deltacheck{\omega}=d\theta$:}
We find, using the transformation rules~\eqref{eq:delf_transf}-
\eqref{eq:m_transf},
\begin{equation}
  \label{eq:deltaomega_tmp}
  \begin{split}
    \deltacheck\omega_{ij} &= \omega_j - \omega_i \\ &=
    \begin{aligned}[t]
    &d\Bigl(2\,\mu_j\,\frac{z^{\prime\prime}_{ij}}{z^\prime_{ij}}
    d\Bar z_j -\bigl(\log (w^\prime_{ij}\circ f_j) + \log
    z^\prime_{ij}\bigr)\,d\log \del f_j + \log
    (w^\prime_{ij}\circ f_j)\, d\log z^\prime_{ij}\Bigr) \\
    &\quad +2\mu_j\,h_j\, dz_j\wedge d\Bar z_j -2
    \mu_i\,h_i\,dz_i\wedge d\Bar z_i
    -2\mu_j\{z_i,z_j\}\,dz_j\wedge d\Bar z_j\,.
    \end{aligned}
  \end{split}
\end{equation}
In light of \eqref{eq:h}, equation \eqref{eq:deltaomega_tmp}
reads
\begin{equation*}
  \deltacheck{\omega}=d\theta,
\end{equation*}
with $\theta$ given by the first two terms on the RHS of
\eqref{eq:deltaomega_tmp}, that is,
\begin{equation}
  \label{eq:theta}
  \begin{split}
    \theta_{ij} =
    2\mu_j\,\frac{z^{\prime\prime}_{ij}}{z^\prime_{ij}} d\Bar z_j
    &- \bigl(\log (w^\prime_{ij}\circ f_j) + \log
    z^\prime_{ij}\bigr)\,d\log \del f_j\\ &+ \log
    (w^\prime_{ij}\circ f_j)\, d\log z^\prime_{ij}.
  \end{split}
\end{equation}

\subsubsection{$\protect\deltacheck\theta$:}
The first term on the RHS of \eqref{eq:theta} is a cocycle, as it
has the \v{C}ech cup product of two terms which are cocycles
themselves.  We can ignore it from now on. The term on the second
line of \eqref{eq:theta} is also cup product, so its coboundary
is computed by applying $\deltacheck(a\cup b)= \deltacheck(a)\cup
b + (-1)^{\mathrm{deg}\,a} a\cup \deltacheck (b)$.  For the
remaining term the cocycle is computed directly. The final result
is
\begin{equation}
  \label{eq:dTheta_1}
  \begin{split}
    \deltacheck\bigl( \theta\bigr)_{ijk} = &-\log w^\prime_{ij}\,
    d\log w^\prime_{jk} +\log z^\prime_{ij}\, d\log
    z^\prime_{jk}\\ &-(\Tilde c_{ijk}+c_{ijk})\,d\log \del f_k
    -d\,\bigl( \log z^\prime_{ij}\log w^\prime_{jk}\bigr)\\
    &\quad +\Tilde c_{ijk}\,d\log z^\prime_{ik}\,,
  \end{split}
\end{equation}
where we suppressed the $f$-dependence. To restore it, notice
that on the triple intersection $U_i\cap U_j\cap U_k$ everything
is evaluated with respect to the coordinate $z_k$, so that $\log
w^\prime_{ij}\circ f_j\vert_{U_k} = \log w^\prime_{ij}\circ
f_j\circ z_{jk} \equiv \log w^\prime_{ij}\circ w_{jk}\circ
f_k$. We shall use this convention in the sequel, in order to
keep some of the expressions less cumbersome.

\subsubsection{$\protect\deltacheck\theta = d\Theta$:}
\label{sec:Theta}
Here we are using Deligne tame symbols in holomorphic category,
introduced in~\ref{sec:prel_sheaves} in order to find $\Theta$
satisfying the equation $\deltacheck\theta = d\Theta$ and to
check that $\deltacheck\Theta=0 \mod\ZZ(2)$.

Consider the tame symbol $\tame{TX}{TX}$, which is represented in
\v{C}ech cohomology by the element
\begin{displaymath}
  \bigl( -\log z^\prime_{ij}\,d\log z^\prime_{jk}\,,\,
  c_{ijk}\log z^\prime_{kl}\,,\, c_{ijk}c_{klm}\bigr) \in
  \Check{C}^2(\sheaf{\Omega}^1_X)\oplus
  \Check{C}^3(\sheaf{\mathcal{O}}_X) \oplus
  \Check{C}^4(\ZZ(2)_X)\,,
\end{displaymath}
where $\{c_{ijk}\}$ represents the first Chern class of $TX$. As
we mentioned in section~\ref{sec:prel_sheaves},
\begin{displaymath}
  \HH^4(X,\ZZ(2)^\bullet_{\mathcal{D},\mathit{hol}}) \cong
  H^3(X,\CC^*)=0,
\end{displaymath}
so that the total cocycle representing $\tame{TX}{TX}$ is a
coboundary:
\begin{displaymath}
  \bigl( -\log z^\prime_{ij}\,d\log z^\prime_{jk}\,,\,
  c_{ijk}\log z^\prime_{kl}\,,\, c_{ijk}c_{klm}\bigr) = D\bigl(
  \tau_{ij}, \phi_{ijk}, n_{ijkl}\bigr)\,,
\end{displaymath}
where $(\tau_{ij})\in
\Check{C}^1(\cover{U}_X,\sheaf{\Omega}^1_X)$, $(\phi_{ijk}) \in
\Check{C}^2(\cover{U}_X,\sheaf{\mathcal{O}}_X)$ and $(n_{ijkl})
\in \Check{C}^3(\cover{U}_X, \ZZ(2)_X)$.  Computing the RHS
yields the relations
\begin{subequations}\label{eq:tame_trivial}
  \begin{align}
    -\log z^\prime_{ij}\,d\log z^\prime_{jk} &= (\deltacheck
     \tau)_{ijk} + d\phi_{ijk}\\ c_{ijk}\log z^\prime_{kl} &=
     -(\deltacheck \phi)_{ijkl} + n_{ijkl}\\ c_{ijk}c_{klm} &=
     (\deltacheck n)_{ijklm}\,.
  \end{align}
\end{subequations}
There is an entirely similar situation for the deformed Riemann
Surface $\Tilde X$ and the corresponding symbol $\tame{T\Tilde
X}{T\Tilde X}$, for which we introduce the corresponding objects
$\Tilde\tau_{ij}$, $\Tilde\phi_{ijk}$ and $\Tilde
n_{ijkl}$. Using these results we rewrite $\deltacheck\theta$ as
\begin{equation*}
  \begin{split}
    \deltacheck\bigl( \theta\bigr)_{ijk}
    &=\deltacheck(f^*(\Tilde\tau))_{ijk} + d
    f^*(\Tilde\phi_{ijk}) -\deltacheck(\tau)_{ijk} -d\phi_{ijk}\\
    &-(\Tilde c_{ijk}+c_{ijk})\,d\log \del f_k -d\,\bigl( \log
    z^\prime_{ij}\log w^\prime_{jk}\bigr)\\ &+\Tilde
    c_{ijk}\,d\log z^\prime_{ik}\,,
  \end{split}
\end{equation*}
where $f^*(\Tilde\tau_{ij})$ and $f^*(\Tilde\phi_{ijk})$ are
pull-backs of forms $\Tilde\tau_{ij}$ and $\Tilde\phi_{ijk}$ on
$X$.  Now, perform the shift:
\begin{displaymath}
  \theta_{ij}\rightarrow \Hat\theta_{ij}\eqdef \theta_{ij}
-f^*(\Tilde\tau_{ij}) +\tau_{ij}\,.
\end{displaymath}
This is possible since $\tau_{ij}$ and $\Tilde\tau_{ij}$ are
holomorphic relative to the respective complex structures,
implying $d\tau_{ij}=0$ and $d f^*(\Tilde\tau_{ij})=0$, so that
\begin{displaymath}
  d\Hat\theta_{ij}=d\theta_{ij}=\deltacheck(\omega)_{ij}\,,
\end{displaymath}
without affecting the 2-form part of the action.

{}From now on we assume that $\theta_{ij}$ has been redefined in
this way, that is
\begin{equation} \label{eq: theta new}
  \theta_{ij} = \theta^{\mathit{old}}_{ij} -f^*(\Tilde\tau_{ij})
  +\tau_{ij}\,,
\end{equation}
where $\theta^\mathit{old}_{ij}$ is given by formula
\eqref{eq:theta}, and we can finally put
\begin{displaymath}
  \deltacheck\theta = d\Theta,
\end{displaymath}
with
\begin{equation}
  \label{eq:Theta}
  \begin{split}
    \Theta_{ijk} &= f^*(\Tilde\phi_{ijk}) - \phi_{ijk} -(\Tilde
c_{ijk}+c_{ijk})\,\log \del f_k\\ &-\log z^\prime_{ij}\,\log
w^\prime_{jk} +\Tilde c_{ijk}\,\log z^\prime_{ik}.
  \end{split}
\end{equation}

\subsubsection{$\protect\deltacheck\Theta$:}
Using the relations \eqref{eq:tame_trivial} we compute:
\begin{equation*}
  \deltacheck\Theta_{ijkl} = \Tilde n_{ijkl} - n_{ijkl} -(\Tilde
c_{ijk} + c_{ijk})b_{kl} + c_{ijl}{\Tilde c}_{jkl} -
c_{ikl}{\Tilde c}_{ijk}\,,
\end{equation*}
so that $\deltacheck\Theta \in\Check{C}^3(\cover{U}_X,\ZZ(2)_X)$.
Setting
\begin{equation*}
  m_{ijkl}\eqdef (\deltacheck\Theta)_{ijkl}\,,
\end{equation*}
we can summarize the foregoing in the following
\begin{proposition}
  \label{prop:tot_cocycle}
  The total cochain
  \begin{equation*}
    \Omega\eqdef 2\pi \sqrt{-1} \bigl( \omega_i,\,
    \theta_{ij},\,
    -\Theta_{ijk},\, -m_{ijkl}\bigr)\,,
  \end{equation*}
  with $\omega_i$ given by the Polyakov form \eqref{eq:omega},
  represents a class in $H^3_\mathcal{D}(X,\ZZ(3))$.
\end{proposition}
\begin{proof} All the preceding computations amount to show
  that
  \begin{equation*}
    D\Omega = 2\pi\sqrt{-1}\bigl( -\omega_j +\omega_i
    +d\theta_{ij},\, (\deltacheck\theta)_{ijk} -d\Theta_{ijk},\,
    (\deltacheck\Theta)_{ijkl} -m_{ijkl},
    (\deltacheck m)_{ijklp}\bigr) = 0\,.
  \end{equation*}
  Then $\Omega$ represents a class since the double complex
  $\CCC^{\bullet,\bullet}_\mathcal{D}$ computes the
  hypercohomology.
\end{proof}

Now that we have constructed the Lagrangian cocycle from the
Polyakov top form in \eqref{eq:omega}, we can finally give the
\begin{definition}
  \label{def:action}
  Let $\mu \in\bbb$ be a Beltrami coefficient, $f$ be the
  associated deformation map, and $\Omega[f]$ be the local
  Lagrangian cocycle constructed from \eqref{eq:omega}. The
  \emph{Polyakov action} functional on $X$ is given by the
  evaluation
  \begin{equation}
    \label{eq:action}
    S[f] \eqdef \langle \Omega [f],\Sigma\rangle\,,
  \end{equation}
  over the representative $\Sigma$ of the fundamental class of
  $X$ given in~\ref{sec:sigma} and in the appendix.
\end{definition}
\begin{remark}
  As it follows from the definition, Polyakov's action is
  well-defined modulo $\ZZ(3)$, so that only its exponential
  $A[f]=\exp\{S[f]/(2\pi i)^2\}$ is well-defined. It also follows
  from the definition of the pairing in section~\ref{sec:sigma}
  that the functional $A[f]$ actually depends only on the
  cohomology class in $H^3_\mathcal{D}(X,\ZZ(3))$ of the local
  Lagrangian cocycle $\Omega[f]$.
\end{remark}
By construction, the cocycle $\Omega[f]$ depends also on a smooth
projective connection $h\in \mathcal{Q}(X)$, so that the
exponential of the action defines the map $A:\mathcal{Q}(X)\times
\bbb\longrightarrow\CC^*$, where the dependence on the first
factor is that of an external field.
\label{sec:properties}

Here we analyze the dependence of the action functionals $S[f]$
and $A[f]$ on the choice of the logarithm branches. We also study
the trivializing coboundary for the tame symbol $\tame{TX}{TX}$,
analyze the dependence of the action on this trivialization, and
show that $A[f]$ should be in fact considered as taking its
values in a $\CC^*$-torsor.

\subsubsection{Dependency on logs}
\label{sec:log_dep}

Here we prove the following

\begin{proposition}
  The functional $A[f]$ is independent of the choice of the
  logarithm branches in \eqref{eq:branches}.
\end{proposition}
\begin{proof}
  It is sufficient to show that changing the definition of the
  various logarithm branches in $\Omega$ amounts to change it by
  a coboundary. First, we change these branches,
  \begin{gather*}
    \log z^\prime_{ij} \longrightarrow \log z^\prime_{ij} +
    k_{ij}\\ \log w^\prime_{ij} \longrightarrow \log
    w^\prime_{ij} + \Tilde k_{ij}\\ \log \del f_i \longrightarrow
    \log \del f_i + p_i
  \end{gather*}
  where $k_{ij},\Tilde k_{ij},p_i\in\ZZ(1)$.  The effect of these
  changes on the representatives of the Chern classes of $TX$ and
  $T\Tilde X$ is
  \begin{gather*}
    b_{ij}\longrightarrow b_{ij} + \Tilde k_{ij} - k_{ij} +p_j
    -p_i\\ c_{ijk}\longrightarrow c_{ijk} +
    \deltacheck(k)_{ijk}\\ \Tilde c_{ijk}\longrightarrow \Tilde
    c_{ijk} +\deltacheck(\Tilde k)_{ijk}\,.
  \end{gather*}
  While the term $\omega_i$ is obviously invariant under these
  changes, $\theta_{ij}$ and $\Theta_{ijk}$, by descent theory,
  transform as follows
  \begin{gather*}
    \theta_{ij}\longrightarrow \theta_{ij} + d\psi_{ij}\,,\\
    \Theta\longrightarrow \Theta +\deltacheck\psi -r_{ijk}\,,
  \end{gather*}
  where $\psi\in\Check{C}^1(\cover{U}_X,\sheaf{A}^0_X)$ and
  $r\in\Check{C}^2 (\cover{U}_X,\CC)$. Note that if $r_{ijk}\in
  \ZZ(2)$ for any $ijk$, then $\Omega\longrightarrow
  \Omega+D\lambda$, where $\lambda= (0,\psi_{ij},r_{ijk})$, and
  we are done.

  To prove that $r\in\Check{C}^2 (\cover{U}_X,\ZZ(2))$, we
  actually compute the shift for $\Theta$. First, we explicitly
  determine
  \begin{equation*}
    \psi_{ij} = -({\Tilde k}_{ij} + k_{ij})\, \log\del f_j
    +\Tilde k_{ij}\, \log z^\prime_{ij}\,.
  \end{equation*}
  Next, we explicitly compute the shift of the total cocycle
  representing $\tame{TX}{TX}$. This is a straightforward
  calculation, using relations \eqref{eq:tame_trivial}, with the
  result:
  \begin{align*}
    \tau_{ij} &\longrightarrow \tau_{ij}\\ \phi_{ijk}
    &\longrightarrow \phi_{ijk} -k_{ij}\,\log z^\prime_{jk}\\
    n_{ijkl} &\longrightarrow n_{ijkl} +k_{ij}c_{jkl}
    +c_{ijk}k_{kl} +(\deltacheck k)_{ijk}k_{kl}\,.
  \end{align*}
  Similar formulas are valid for the shift of $\tame{T\Tilde
  X}{T\tilde X}$. Putting everything together, we get
  \begin{multline}
    \label{eq:Theta_shift_1}
    r_{ijk}= (\Tilde k_{ij} +k_{ij})b_{jk} +\bigl(\Tilde c_{ijk}
    +c_{ijk} +(\deltacheck k)_{ijk} +(\deltacheck\Tilde
    k)_{ijk}\bigr)p_k\\ +k_{ij}\Tilde k_{jk} -\Tilde
    c_{ijk}k_{ik} -(\deltacheck \Tilde k)_{ijk}k_{ik} +\Tilde
    k_{jk}c_{ijk} +\Tilde k_{ij}c_{ijk}\in\ZZ(2).
  \end{multline}
\end{proof}

\subsubsection{A more detailed analysis of the vanishing tame symbol}
\label{sec:bloch_dilog}

Here we analyze the condition $\tame{TX}{TX}=0$ as an element of
$\HH^4(X,\ZZ(2)^\bullet_{\mathcal{D},\mathit{hol}})$ in more
detail. In particular, we investigate the possibility of putting
the trivializing cochains $(\tau_{ij},\phi_{ijk},n_{ijkl})$ and
$(\Tilde\tau_{ij},\Tilde\phi_{ijk},\Tilde n_{ijkl})$ into some
specific forms. This analysis is based on the
relations~\eqref{eq:tame_trivial}, which we rewrite here:
\begin{align*}
    -\log z^\prime_{ij}\,d\log z^\prime_{jk} &= (\deltacheck
     \tau)_{ijk} + d\phi_{ijk}\\ c_{ijk}\log z^\prime_{kl} &=
     -(\deltacheck \phi)_{ijkl} + n_{ijkl}\\ c_{ijk}c_{klm} &=
     (\deltacheck n)_{ijklm}\,.
\end{align*}

The first equation above calls for the differential equation
\begin{displaymath}
    -\log z^\prime_{ij}\circ z_{jk}\,d\log z^\prime_{jk} =
     dL_{ijk}\,.
\end{displaymath}
Its solution $L_{ijk}(z_k)$ can be considered as a Bloch
dilogarithm associated to the symbol
$\tame{z^\prime_{ij}}{z^\prime_{jk}}$, which is the cup-product
in Deligne cohomology of the two invertible functions
$z^\prime_{ij}$ and $z^\prime_{jk}$ and is a trivial element of
$H^2_\mathcal{D}(U_{ijk},\ZZ(2))$ (see~\cite{esn-vie} for more
details).  The consistency condition on quadruple intersections
$U_{ijkl}$ is obtained by applying the \v{C}ech coboundary to the
differential equation satisfied by $L_{ijk}$. One gets
\begin{displaymath}
     c_{ijk}\log z^\prime_{kl} = -(\deltacheck L)_{ijkl} +
     \alpha_{ijkl}\,,
\end{displaymath}
where $\alpha_{ijkl}$ is a $\CC$-valued cochain --- an
integration constant. By taking the \v{C}ech coboundary of the
last relation we get
\begin{displaymath}
     c_{ijk}c_{klm} = (\deltacheck \alpha)_{ijklm}\,.
\end{displaymath}
Therefore,
\begin{displaymath}
  \deltacheck (\alpha - n)=0\,,
\end{displaymath}
that is, the element $\alpha -n$ is a 3-cocycle. By dimensional
reasons, it must be a coboundary,
\begin{displaymath}
  \alpha = n+\deltacheck\beta\,,
\end{displaymath}
with $\beta$ being a $2$-cochain with values in $\CC$. It follows
that
\begin{displaymath}
  c_{ijk}\log z^\prime_{kl} = -\deltacheck( L-\beta)_{ijkl} +
  n_{ijkl}\,.
\end{displaymath}

As a result, we effectively obtained a trivializing cocycle for
the tame symbol $(TX,TX]$ which does not include a 1-form:
\begin{displaymath}
  \bigl( -\log z^\prime_{ij}\,d\log z^\prime_{jk}\,,\,
  c_{ijk}\log z^\prime_{kl}\,,\, c_{ijk}c_{klm}\bigr) =
  D\bigl(0,L_{ijk},n_{ijkl} \bigr)\,,
\end{displaymath}
where we relabeled $L-\beta\rightarrow L$.

\subsubsection{Relation with $\CC^*$-torsors}
\label{sec:torsor}

Notice that the trivialization of the tame symbol $\tame{TX}{TX}$
is defined up to a cocycle representing an element in
$\HH^3(X,\ZZ(2)^\bullet_{\mathcal{D},\mathit{hol}}) \cong
H^2(X,\CC/\ZZ(2)) \cong H^2(X,\CC^*) \cong \CC^*$.  Thus there
is a $\CC^*$-action on the functional $A[f]$ which simply is the
shift of the total trivializing cochain
$(\tau_{ij},\phi_{ijk},n_{ijkl})$ by a cocycle representing a
class in $\HH^3(X,\ZZ(2)^\bullet_{\mathcal{D},\mathit{hol}})$.
From this it is clear that, keeping $f$ fixed, the functional
$A[f]$ does not simply take its values in $\CC^*$, but rather in
a $\CC^*$-torsor $T$. From this perspective, choosing a specific
total cochain to trivialize the symbol $\tame{TX}{TX}$ amounts to
choosing an isomorphism $T\overset{\sim}{\rightarrow} \CC$.

The $\CC^*$-action can be described explicitly if we make use of
the cocycle $(0,L_{ijk},n_{ijkl})$, obtained by choosing a
dilogarithm $L_{ijk}$ for the symbol
$\tame{z^\prime_{ij}}{z^\prime_{jk}}$. Namely, as it follows from
the discussion in the previous section, we can add to $L_{ijk}$ a
cocycle $(\beta_{ijk},p_{ijkl})$ representing an element in
\begin{displaymath}
  \HH^3(X,\ZZ(2)\overset{\imath}{\rightarrow} \CC)\cong
  H^2(X,\CC^*)\,.
\end{displaymath}
Note that, by definition, $\deltacheck\beta = p\in \ZZ(2)$.

Since the action functional is defined using trivialization of
\emph{two} tame symbols, $\tame{TX}{TX}$ and $\tame{T\Tilde
  X}{T\Tilde X}$, the above argument should be applied to both
cochains $(\tau_{ij},\phi_{ijk},n_{ijkl})$ and
$(\Tilde\tau_{ij},\Tilde\phi_{ijk},\Tilde n_{ijkl})$, so that we
have in fact two $\CC^*$-actions. From a Teichm\"uller theory
point of view, these two actions refer to very different
structures. One is defined in terms of the complex structure $X$
which is fixed throughout (a base point in Teichm\"uller space),
while the other is relative to the $f$-dependent complex
structure $\Tilde X$. The latter action depends on the dynamical
field $f$.

Thus it is appropriate to speak of a $(\CC^*,\CC^*)$-action, in
the sense that the space $T$ where the action takes its values
carries two simultaneous (and compatible) $\CC^*$-actions.

\subsection{Other coverings --- a dictionary}
\label{sec:dictionary}
In this section we set up a dictionary connecting generalized \v
Cech formalism developed in \ref{sec:prel_coverings} and
\ref{sec:prel_fundamental} with the formalism used
in~\cite{prev_paper} for the universal cover of $X$. Besides
comparing the two formalisms, by applying the dictionary to the
formulas in \ref{sec:local_cocycle}, we also clarify the explicit
form of the Lagrangian cocycle constructed in~\cite{prev_paper}.
Specifically, we treat the ``integration constants'' arising from
solving the descent equations via Deligne complexes and analyze
explicit dependence of the action functional on background
projective structures.

\subsubsection{}
Start from the universal cover $U\to X$, which we specify as the
upper half-plane $\HH$. Then $\mathrm{Deck}(\HH/X)\cong
\pi_1(X)\cong \Gamma$, a finitely-generated, purely hyperbolic
Fuchsian group (a discrete subgroup of $\PSL_2(\RR)$),
uniformizing the Riemann surface $X$. The group $\Gamma$ acts on
$\HH$ by M\"obius transformations.

Geometric objects on $X$ correspond to $\Gamma$-equivariant
objects on $\HH$: a tensor $\phi\in A^{p,q}(X)$ corresponds to an
automorphic form $\phi$ for $\Gamma$ of weight $(2p,2q)$, i.e.~a
function (indicated by the same name) $\phi:\HH\to \CC$ such that
\begin{equation*}
\phi\circ\gamma \cdot (\gamma^\prime )^p
(\overline{\gamma^\prime})^q = \phi,\,\gamma\in\Gamma.
\end{equation*}
Clearly, an automorphic form is just a zero cocycle on $\Gamma$
with values in $\sheaf{A}^{p,q}(\HH)$. Examples of automorphic
forms of geometric origin are provided by Beltrami differentials
on $X$, that correspond to forms of weight $(-2,2)$, by abelian
differentials on $X$ --- global sections of $\sheaf{\Omega}_X$
--- that correspond to holomorphic forms of weight $(2,0)$, and
by quadratic differentials on $X$ --- global sections of
$\sheaf{\Omega}^{\otimes2}_X$ --- that correspond to holomorphic
forms of weight $(4,0)$.

The deformation map $f$ is realized as a quasi-conformal map
\begin{equation*}
  f:\HH \longrightarrow \DD
\end{equation*}
satisfying on $\HH$ the Beltrami equation
\begin{equation}\label{eq:beltrami}
  f_{\bar z} = \mu f_z\,,
\end{equation}
where $\mu$ is a Beltrami differential for $\Gamma$ on $\HH$ such
that $\Vert\mu\Vert_{\infty}<1$. The Beltrami equation on $\HH$
should be supplemented by boundary conditions, that guarantee the
following.
\begin{enumerate}
\item $\DD=f(\HH)$ is a \emph{quasi-disk}, i.e.~a domain in
  $\PP^1$ bounded by a closed Jordan curve and analytically
  isomorphic to $\HH$;
\item $\Tilde\Gamma = f\circ \Gamma\circ f^{-1}\subset
  \PSL_2(\CC)$ is a discrete subgroup, isomorphic to $\Gamma$ as
  an abstract group, acting on $\DD$, i.e.~a so-called
  \emph{quasi-fuchsian} group. The isomorphism $\Gamma\to
  \Tilde\Gamma$ intertwines $f$.
\end{enumerate}
These boundary conditions are specified by extending $\mu$ to the
whole complex plane, where the Beltrami equation has a unique
solution up to a post-composition with M\"obius
transformation~\cite{ahl,nag}. The following two types are of
particular importance.
\begin{enumerate}
\item[(a)] Extension of $\mu$ by reflection to the lower half
  plane $\Bar\HH$: $\mu (z)\eqdef \overline{\mu(\Bar z)}$ for
  $z\in \overline\HH$. Then $\DD=\HH$ and $\Tilde\Gamma$ is also
  a Fuchsian group.
\item[(b)] Extension of $\mu$ by setting $\mu (z) = 0$ for $z\in
  \overline\HH$. In this case $\DD$ is a quasi-disc and the
  dependence of the mapping $f$ on $\mu$ is holomorphic.
\end{enumerate}
The formalism developed below will be independent of a particular
boundary condition chosen.

\subsubsection{}
Here we address a minor normalization problem caused by the fact
that the action of $\PSL_2(\RR)$ --- and therefore of $\Gamma$
and $\Tilde\Gamma$ --- by M\"obius transformations is on the left
instead of on the right, as we assumed in
\ref{sec:prel_coverings}. Assuming a right action yields all the
standard formulas in group cohomology. On the other hand, a left
action of $\Gamma$ is more convenient in view of the fact that
$\HH$ itself is the quotient of a principal fibration:
$\HH\cong\PSL_2(\RR)/\SO(2)$.\footnote{In general, we prefer to
consider \emph{right} principal fibrations.} As a result, the
surface itself is presented as a double coset space: $X\cong
\Gamma\backslash \PSL_2(\RR)/\SO(2)$.

For a left action $G\times U\to U$ for a $G$-space $U\to M$ there
is the isomorphism
\begin{equation*}
  \underbrace{U\times_M\dots\times_M U}_{q+1} \cong G^q \times
  U\,,
\end{equation*}
sending the $q$-tuple $(x_0,\dots,x_q)$ to the tuple
$(g_1,\dots,g_q,x)$ such that
\begin{equation*}
  (x_0,\dots,x_q) = (g_1\dots g_q x,g_2\dots g_q x,\dots, g_q x,
  x).
\end{equation*}
This arrangement makes the face maps $d_i$ appear in backward
order, that is
\begin{equation*}
  d_0(g_1,\dots,g_q,x)=(g_2,\dots,g_q,x) \,\dots\,
  d_q(g_1,\dots,g_q,x)=(g_1,\dots,g_{q-1},g_qx)\,.
\end{equation*}
As a result, the action on the coefficients would be on the right
and the coboundary operator $\deltapp$ in group cohomology should
actually be read from right to left, as in
\begin{equation}\label{eq:non_std_coboundary}
  \begin{split}
    (\deltapp\phi)_{g_1,\dots,g_q} = \phi_{g_2,\dots,g_q}
    &+\sum_{i=1}^{q-1}(-1)^i
    \phi_{g_1,\dots,g_ig_{i+1},\dots,g_q}\\
    &+(-1)^q{g_q}^*\phi_{g_1,\dots,g_{q-1}}\,,
  \end{split}
\end{equation}
for $\phi$ a $(q-1)$ cochain. Observe that the pull-back action
on the coefficients is a right one.

The familiar formulas in group cohomology can be retrieved by
turning the left action into a right one using the standard trick
\begin{equation*}
  x\cdot g \eqdef g^{-1}x\,,\quad g\in G,x\in U\,,
\end{equation*}
which at the level of nerves amounts to performing the swap
$(g_1,\dots,g_q,x) \mapsto (x, g_q^{-1},\dots, g_1^{-1})$ in
degree $q$. It follows that one has to evaluate all cochains over
inverses of group elements. This is the convention we followed
in~\cite{prev_paper}.

On the other hand, given the action of $\Gamma$ on $\HH$ as a
left one, keeping the non standard form
\eqref{eq:non_std_coboundary} parallels more closely the \v Cech
framework if we consider the pair $(\gamma(z),z)\in
\HH\times_X\HH$, for $z\in \HH$ and $\gamma\in \Gamma$, as a
change of coordinates, much like the pair $(z_i,z_j)\in
U_i\times_X U_j\equiv U_i\cap U_j$ with $z_i = z_{ij}(z_j)$.
More generally, for $U_{i_0}\cap\dots\cap U_{i_q}$ there are
coordinates $z_{i_0},\dots ,z_{i_q}$ with $z_{i_k}=z_{i_k
i_{k+1}}(z_{i_{k+1}})$, $k=0,\dots,q-1$, and if $\phi
\in\Check{C}^{q-1}(\cover{U}_X,\sheaf{A}^p)$ then we have
\begin{equation}\label{eq:cech_pedantic}
  \deltacheck\phi_{i_0,\dots,i_q} = \sum_{k=0}^{q-1} (-1)^k
  \phi_{i_0,\dots,\Hat{i_k},\dots,i_q} +(-1)^q
  (z_{i_{q-1}i_q})^*\phi_{i_0,\dots,i_{q-1}}\,,
\end{equation}
where the convention is that each component 
is expressed in the coordinate 
determined by the last index.  This is the formula we used when
performing explicit computations with \v Cech cochains for the
calculation of the local Lagrangian cocycle. Thus
\eqref{eq:cech_pedantic} becomes formally equal to
\eqref{eq:non_std_coboundary} when we interpret the last
pull-back by $g_q$ as the restriction isomorphism expressing
everything in terms of the last coordinate.

\subsubsection{}
The translation of the constructions in \ref{sec:action_cech} and
\ref{sec:local_cocycle} to the upper-half plane is now done
according to the following table:
\begin{center}
  \begin{tabular}{|c|c|}
    \hline \v Cech:$\cover{U}_X$ & Upper-half plane $\HH$\\
    \hline \hline $U_{i_0}\cap\dots\cap U_{i_n}$ &
    $\Gamma^n\times \HH$\\ \hline $z_{i_0},\dots ,z_{i_n}$ &
    $\gamma_1,\dots,\gamma_n,z$\\ \hline $z_{i_k}=z_{i_k
    i_{k+1}}(z_{i_{k+1}}),k=0,\dots,n-1$ &
    $z_{k-1}=\gamma_k(z_k),k=1,\dots,n,z_n=z$\\ \hline
    $\phi_{i_0,\dots, i_n}(z_{i_n})dz_{i_n}^p d\Bar z_{i_n}^q$ &
    $\phi_{\gamma_1,\dots,\gamma_n}(z)dz^p d\Bar z^q$\\
    \hline\eqref{eq:cech_pedantic} &
    \eqref{eq:non_std_coboundary}\\ \hline
  \end{tabular}
\end{center}
Similar provisions of course relate the deformed coordinates
$w_i$ and elements of the deformed group $\Tilde\Gamma$. Any
construction explicitly involving the map $f$ must take into
account the equivariance property $f\circ\gamma = \Tilde\gamma
\circ f$ for any $\gamma\in \Gamma$, where $\Tilde\gamma$ is the
corresponding element in the deformed group $\Tilde\Gamma$. We
have relations entirely similar to \eqref{eq:delf_transf} and
\eqref{eq:delbf_transf} which can be found in~\cite{prev_paper};
for example
\begin{equation}\label{eq:delf_transf.3}
  \frac{\Tilde\gamma^\prime\circ f}{\gamma^\prime} f_z = f_z\circ
  \gamma\,.
\end{equation}
In order to handle the logarithm of \eqref{eq:delf_transf.3} in
the same way as we just did in the \v Cech case (see
\eqref{eq:branches}) we depart from~\cite{prev_paper}. The
problem is to relate $\log (\gamma_1\gamma_2)^\prime$ and
$\log\gamma_1^\prime\circ \gamma_2 +\log\gamma_2^\prime$ for any
$\gamma_1, \gamma_2\in \Gamma$, and similarly for $\Tilde\Gamma$.
Instead of directly analyzing the branch-cuts (thus introducing
an element of explicit dependence on the choice of the branches)
we set
\begin{subequations}
  \label{eq:branches.2}
  \begin{align}
    c_{\gamma_1,\gamma_2} &= \log \gamma_2^\prime -\log
(\gamma_1\gamma_2)^\prime +\log\gamma_1^\prime\circ\gamma_2\,,\\
\Tilde c_{\Tilde\gamma_1,\Tilde\gamma_2} &= \log
\Tilde\gamma_2^\prime - \log(\Tilde\gamma_1\Tilde\gamma_2)^\prime
+\log\Tilde\gamma_1^\prime\circ \Tilde\gamma_2\,,\\ b_{\gamma} &=
\log \Tilde\gamma^\prime\circ f - \log \gamma^\prime -\log
f_z\circ \gamma +\log f_z \,.
  \end{align}
\end{subequations}
The numbers $c_{\gamma_1,\gamma_2}$, $\Tilde
c_{\Tilde\gamma_1,\Tilde\gamma_2}$ and $b_\gamma$ belong to
$\ZZ(1)$, and $c$, $\Tilde c$ are cocycles with $\Tilde c= c
+\deltapp b$. Since $\gamma^\prime$ is the automorphy factor for
$TX$, the geometric interpretation is that again $c$ represents
$c_1(TX)$~\cite{gun-theta}. Alternatively, $c$ represents the
Euler class of the $S^1$-bundle $\Gamma \backslash \PSL_2(\RR)\to
X$~(\cite{milnor,wood,goldman}, see
also~\cite{matsumoto}). Indeed, the first of equations
\eqref{eq:branches.2} can be written in terms of rotation
numbers:
\begin{equation*}
  c_{\gamma_1,\gamma_2}= -2\bigl(w(\gamma_2) -w(\gamma_1\gamma_2)
+w(\gamma_1)\circ\gamma_2\bigr)\,,
\end{equation*}
where $w\bigl(
\begin{smallmatrix}
  a & b\\ c&d
\end{smallmatrix}\bigr) = \arg (c z +d)$. More precisely, this is
the Euler class of the $\RR\PP^1$-bundle obtained by letting
$\PSL_2(\RR)$ act on the real projective line realized as the
boundary of $\HH$ (see~\cite{wood} for details). Again, a similar
discussion holds for $\Tilde\Gamma$ with the obvious changes.

As was shown in section~\ref{sec:prel_coverings}, \v Cech
cohomology with respect to the cover $\HH\to X$ is the same as
group cohomology of $\pi_1(X)\cong \Gamma$ with values in the
appropriate coefficients. Also it was noted there that $\HH\to X$
is a good covering acyclic for fine sheaves, so that
$H^p(\pi_1(X),\CC)\cong H^p_\mathit{dR}(X)\cong H^p(X,\CC)$.
Similar arguments show that the double complex
$C^q(\Gamma,\ZZ(3)^p_\mathcal{D})$ computes
$H^\bullet_\mathcal{D}(X,\ZZ(3))$.

The choice of the covering $\HH\to X$ --- or, more generally,
$\DD\to X$ --- contains more information than simply using an
abstract universal covering map $U\to X$: it includes the choice
of a projective structure. Indeed, since the Schwarzian
derivative of any M\"obius transformation vanishes, any local
section of the canonical projection would precisely be a system
of projective charts for it.

It follows that when working with $\HH\to X$ the explicit
inclusion of projective connections becomes --- strictly speaking
--- unnecessary. Indeed, these were not considered
in~\cite{prev_paper}. However, it is
known~\cite{prev_paper,zucchini} that the effective action (that
is, the class of the local Lagrangian cocycle) in higher genus is
determined --- say, by the Universal Ward Identity --- only up to
holomorphic quadratic differentials. Interpreting the latter as
lifts of projective connections, the precise statement is that
the effective action is determined up to the choice of a
projective structure.  In light of this observation, and also to
keep a strict parallel with the \v Cech formulation, we make this
dependence on a generic projective connection
explicit.\footnote{Here the term projective connection is to be
  understood in the same way as in
  \protect\ref{sec:local_cocycle}, i.e.~as not necessarily
  holomorphic one.}  In this way we obtain a unified formalism
consistent with the treatment of variations in section
\ref{sec:variation_proj}, where conditions on the projective
connections will be enforced by the variational process.

Now we set out to write the correspondence:
\begin{center}
  \begin{tabular}{|c||c|c|}
    \hline & $\cover{U}_X$ & $\HH$\\ \hline\hline $(3,0)$ &
    $\omega_i(z_i)$ & $\omega (z)$ \\ \hline $(2,1)$ &
    $\theta_{ij}(z_j)$ & $\theta_{\gamma}(z)$\\ \hline $(1,2)$ &
    $\Theta_{ijk}(z_k)$ & $\Theta_{\gamma_1,\gamma_2}(z)$ \\
    \hline $(0,3)$ & $m_{ijkl}$ & $m_{\gamma_1, \gamma_2,
    \gamma_3}$ \\ \hline
  \end{tabular}
\end{center}
For the first two lines we start by translating \eqref{eq:omega}
and \eqref{eq:theta}, respectively:
\begin{gather}
  \label{eq:transl_omega}
  \omega_\gamma(z)=\frac{\del^2f}{\del f}\, \del\mu\, dz\wedge
  d\Bar z +2\,\mu\,h\,dz \wedge d\Bar z\,, \\
  \label{eq:transl_theta_first}
  \theta_{\gamma}(z) =
  2\mu\,\frac{\gamma^{\prime\prime}}{\gamma^\prime} d\Bar z -
  \bigl(\log ({\Tilde\gamma}^\prime\circ f) + \log
  \gamma^\prime\bigr)\,d\log \del f + \log
  ({\Tilde\gamma}^\prime\circ f)\, d\log \gamma^\prime\,,
\end{gather}
where $h$ is a smooth quadratic differential. In this way the
last term of \eqref{eq:transl_omega} is automorphic of weight
$(1,1)$, hence it is killed by the coboundary operator. This
would be consistent with a translation of
\eqref{eq:deltaomega_tmp}.  We stress
\eqref{eq:transl_theta_first} is a direct translation of the
expression for the $(2,1)$ component \emph{prior} to the
computation of $\deltapp\theta = d\Theta$. As before, the
existence of $\Theta_{\gamma_1,\gamma_2}$ is guaranteed by the
vanishing of the analog of the symbol $\tame{TX}{TX}$ in
holomorphic Deligne cohomology. This time, the tame symbol is
represented by the cocycle
\begin{equation*}
    \bigl( -\log \gamma_1^\prime\circ \gamma_2 \,d\log
  \gamma_2^\prime\,,\, c_{\gamma_1,\gamma_2}\log
  \gamma_3^\prime\,,\, c_{\gamma_1,\gamma_2}
  c_{\gamma_3,\gamma_4}\bigr)\in
  C^2(\Gamma,\sheaf{\Omega}^1(\HH))\oplus
  C^3(\Gamma,\sheaf{\mathcal{O}}(\HH))\oplus
  C^4(\Gamma,\ZZ(2))\,.
\end{equation*}
Since $\HH\to X$ is a good cover, the quasi-isomorphism
\begin{displaymath}
  \ZZ(2)^\bullet_{\mathcal{D},\mathit{hol}}
  \overset{\sim}{\rightarrow} \bigl( \sheaf{\mathcal{O}}^*(\HH)
  \xrightarrow{d\log} \sheaf{\Omega}^1(\HH) \bigr)[-1]
  \overset{\sim}{\rightarrow} \CC/\ZZ(2) \cong \CC^*
\end{displaymath}
is still in place by holomorphic Poincar\'e lemma on $\HH$.
Hence
\begin{displaymath}
  \HH^4(\Gamma,\ZZ(2)^\bullet_{\mathcal{D},\mathit{hol}}) \cong
  H^3(\Gamma,\CC^*)=0 \,,
\end{displaymath}
again, by obvious dimensional reasons. It follows that we can
still introduce $(\tau_{\gamma})\in
C^1(\Gamma,\sheaf{\Omega}^1(\HH))$, $(\phi_{\gamma_1,\gamma_2})
\in C^2(\Gamma,\sheaf{\mathcal{O}}(\HH))$ and
$(n_{\gamma_1,\gamma_2,\gamma_2}) \in C^3(\Gamma, \ZZ(2))$ such
that
\begin{equation*}
  \bigl( -\log \gamma_1^\prime\circ \gamma_2 \,d\log
  \gamma_2^\prime\,,\, c_{\gamma_1,\gamma_2}\log
  \gamma_3^\prime\,,\, c_{\gamma_1,\gamma_2}
  c_{\gamma_3,\gamma_4}\bigr) = D\bigl( \tau_{\gamma},
  \phi_{\gamma_1,\gamma_2},
  n_{\gamma_1,\gamma_2,\gamma_3}\bigr)\,.
\end{equation*}
where various $\gamma_i$'s are used as place-holders for added
clarity. Obviously, the treatment for the corresponding
quantities depending on $\Tilde\Gamma$ is entirely similar. As a
result, we can either compute the coboundary of
\eqref{eq:transl_theta_first} or simply translate
\eqref{eq:dTheta_1} and repeat step by step what was done in
section \ref{sec:local_cocycle} to arrive at
\begin{equation}
  \label{eq:transl_theta}
  \theta_{\gamma} = \theta^{\mathit{old}}_{\gamma}
-\Tilde\tau_{\gamma} +\tau_{\gamma}\,,
\end{equation}
with $\theta^{\mathit{old}}_{\gamma}$ given by
\eqref{eq:transl_theta_first} and, finally:
\begin{gather}
  \label{eq:transl_Theta}
    \begin{split}
      \Theta_{\gamma_1,\gamma_2} &=
      \Tilde\phi_{\gamma_1,\gamma_2} -\phi_{\gamma_1,\gamma_2}
      -(\Tilde c_{\gamma_1,\gamma_2}+c_{\gamma_1,\gamma_2})\,
      \log\del f\\ &-\log\gamma_1^\prime\circ\gamma_2\, \log
      {\Tilde\gamma_2}^\prime +\Tilde c_{\gamma_1,\gamma_2}\,
      \log (\gamma_1\circ\gamma_2)^\prime
    \end{split}\\
    \label{eq:transl_m}
    \begin{split}
      m_{\gamma_1,\gamma_2,\gamma_3} &= \Tilde
      n_{\gamma_1,\gamma_2,\gamma_3} -
      n_{\gamma_1,\gamma_2,\gamma_3} -(\Tilde
      c_{\gamma_1,\gamma_2}
      +c_{\gamma_1,\gamma_2})\,b_{\gamma_3}\\
      &+c_{\gamma_1,\gamma_2\circ\gamma_3} {\Tilde
      c}_{\gamma_2,\gamma_3} -c_{\gamma_1\circ\gamma_2,\gamma_3}
      {\Tilde c}_{\gamma_1,\gamma_2}\,.
  \end{split}
\end{gather}
Therefore the analog of proposition \ref{prop:tot_cocycle} holds
\begin{proposition}
  \label{prop:transl_tot_cocycle}
  The total cochain
  \begin{equation*}
    \Omega\eqdef 2\pi\sqrt{-1}\bigl( \omega,\, \theta_{\gamma},\,
-\Theta_{\gamma_1,\gamma_2},\,
-m_{\gamma_1,\gamma_2,\gamma_3}\bigr)\,,
  \end{equation*}
  with $\omega$ given by the Polyakov form
  \eqref{eq:transl_omega}, represents a class in
  $H^3_\mathcal{D}(X,\ZZ(3))$.
\end{proposition}

The action functional $S[f]$ is computed by evaluating
$\Omega[f]$ over the appropriate representative of $[X]$, which
in this case would be a total cocycle in
$S_p(\HH)\otimes_{\ZZ\Gamma}B_q(\Gamma)$ whose $(2,0)$ component
can be taken as a fundamental domain $F$ for $\Gamma$ in the form
of a standard $4g$-gon, as detailed in~\cite{prev_paper}.

\section{Variation and projective structures}
\label{sec:variation_proj}

\subsection{Variation}
\label{sec:variation}
Here we compute the variation of the action functional $S[f]$
with respect to the dynamical field $f$, i.e.~we compute its
differential in field space. We denote by $\var$ the variational
operator --- the exterior differential in field space~
\cite{takens,zuckerman,del-freed} --- and we will use coordinates
with respect to a good \v{C}ech cover $\cover{U}_X$ whenever a
local computation is required.

Since the dynamical field $f$ is a deformation map on $X$, we can
either choose to allow variations that effectively deform the
complex structure or restrict ourselves to the ``trivial'' ones
--- deformations corresponding to vertical tangent vectors in the
Earle-Eells fibration over the Teichm\"uller space.

{}From \eqref{eq:f_transf} we get
\begin{equation}
  \label{eq:KS}
  f^*(\kappa_{ij}) = \frac{\var f_i}{\del f_i}\circ z_{ij}\cdot
  (z^\prime_{ij})^{-1} - \frac{\var f_j}{\del f_j}\,,
\end{equation}
where
\begin{displaymath}
\kappa=\{\kappa_{ij}\eqdef \frac{\var w_{ij}}{w^\prime_{ij}}\}
\end{displaymath}
is the standard Kodaira-Spencer deformation cocycle, and
$f^*(\kappa_{ij})=\kappa_{ij}\circ f_j/\del f_j$ is its
pull-back.  The condition $[\kappa]=0$ in $H^1(\Tilde
X,\Tilde\Theta)$, where $\Tilde\Theta$ is the tangent sheaf of
$\Tilde X$, selects variations that leave the complex structure
$X$ fixed. Specifically, if $[\kappa]=0$ then it follows from
\eqref{eq:KS} that $\var f_i/\del f_i$ represents a smooth
$(1,0)$-vector field on $X$ --- possibly after redefining it by a
holomorphic coboundary for $f^*(\kappa_{ij})$. Furthermore, the
variation $\var\mu$ of the corresponding Beltrami differential as
a tangent vector to $\bbb$ at $\mu$ is
\begin{equation}
  \label{eq:var_mu}
  \var\mu = \delb_\mu \frac{\var f}{\del f}\,,
\end{equation}
so the class $[\var\mu]\in H^{(-1,1)}_{\delb_\mu}(X)$ corresponds
to $[\kappa]$ under the Dolbeault isomorphism.

In the sequel we shall confine ourselves to \emph{vertical
  variations,} that is, to those with $[\kappa]=0$.  Then
$\frac{\var f_i}{\del f_i}$ defines a smooth vector field on $X$.

We start to compute the variation of the Lagrangian cocycle
$\Omega$ with respect to $f$.  From a purely formal point of
view, the calculation for the variation of the top form part
proceeds as usual, where in each coordinate patch we have
\begin{equation*}
  \var\omega_i = a_i + d\eta_i\,
\end{equation*}
with $(a_i)\in \Check{C}^0(\cover{U}_X,\sheaf{A}^2_X)$ and
$(\eta_i)\in \Check{C}^0(\cover{U}_X,\sheaf{A}^1_X)$, where
\begin{equation*}
  a_i(f,\var f) = -2\,\delb_\mu \bigl(%
  h_i -\{f_i,z_i\}\bigr)\,\frac{\var f_i}{\del f_i}\,
  dz_i \wedge d\Bar z_i\,.
\end{equation*}
Using the well-known identity
\begin{equation*}
  \delb_\mu \{f,z\} = \del^3\mu\,,
\end{equation*}
where $\mu=\mu(f)$ and $z$ is a local coordinate on $X$ (the
index $i$ is omitted here), we get
\begin{align*}
  \delb_\mu(\{f,z\} -h) &= (\delb -\mu\,\del -2\,\del\mu)(\{f,z\}
  -h) \\ &= \del^3\mu - (\delb -\mu\,\del -2\,\del\mu)h \\ &=
  \del^3\mu +2h\,\del \mu + \del h\,\mu - \delb h\\ &=
  \diff{D}_h\mu - \delb h\,.
\end{align*}
Here, for any smooth projective connection $h\in\mathcal{Q}(X)$,
$\diff{D}_h$ is the following third order differential operator:
\begin{equation}
  \label{eq:delcube}
  \diff{D}_h=\del^3 +2\,h\,\del +\del h\,.
\end{equation}
It is well-known (see, e.g.,~\cite{gunning-paper}) that it has
the property
\begin{equation*}
  \diff{D}_h : \sheaf{A}^{-1,l}_X \longrightarrow
  \sheaf{A}^{2,l}_X\,,
\end{equation*}
for all $l$; in particular, $\diff{D}_h$ maps global forms of
weight $(-1,l)$ to global forms of weight $(2,l)$.

Thus the final expression for the variation of the top form term
is,
\begin{equation}
  \label{eq:var_a}
  a_i(f,\var f) = -2\, \bigl( \delb h_i -\diff{D}_h\,\mu_i
  \bigr)\,\frac{\var f_i}{\del f_i} \, dz_i \wedge d\Bar z_i\,.
\end{equation}
Thanks to \eqref{eq:h} and to the fact that $\diff{D}_h$ is a
well defined map, $a_i(f,\var f)$ is a well defined global
$2$-form on $X$.  The 1-form $\eta_i$ has the expression
\begin{equation}
  \label{eq:var_eta}
  \eta_i = \var\log \del f_i\,d\log\del f_i +2 \del(\log \del
  f_i)dz_i\,\interior\, \var\mu_i -2\,(h_i
  -\{f_i,z_i\})\,\frac{\var f_i}{\del f_i}\, \bigl(dz_i +\mu_i
  d\Bar z_i\bigr)\,,
\end{equation}
where $\interior$ is the interior product between 1-forms and
vectors.

The main point is that the term \eqref{eq:var_a} alone
constitutes the variation of the \emph{whole} Lagrangian cocycle.
Namely, we have
\begin{theorem}
  \label{thm:variation}
  The variation of the total cocycle $\Omega[f] = 2\pi\sqrt{-1}
  \bigl( \omega_i,\,\theta_{ij},\, -\Theta_{ijk},\,
  -m_{ijkl}\bigr)$ under vertical variation is given by the
  2-form \eqref{eq:var_a} up to a total coboundary in the Deligne
  complex. The variation of the action functional $S[f]$ is
  \begin{equation*}
    \var S[f] = 2\pi\sqrt{-1} \int_X a(f,\var f)\,,
  \end{equation*}
  giving the following Euler-Lagrange equation
  \begin{equation}
    \label{eq:proj_condition}
    \diff{D}_h\mu - \delb h=0\,.
  \end{equation}
\end{theorem}
We can give two different proofs of this theorem. One is more in
keeping with the spirit of this work and uses the explicit form
of $\Omega$. The other is based only on Takens' acyclicity
theorem~\cite{takens} for the variational bicomplex and the
formal machinery of descent equations. Although we present both,
the second one will only be sketched here, as providing details
for it would lead us to far afield.\footnote{We plan to return to
  the topic from a more general point of view elsewhere.}
\begin{proof}[First proof]
  The procedure is to compute the variation of the various
  components of $\Omega$ by applying $\var$ to the descent
  equations. Start with $\deltacheck\var\omega_{ij}$, that can be
  computed in two different ways: from equation
  $\deltacheck\omega=d\theta$, and from the variational relation
  $\var \omega = a +d\eta$. Since $a_i=a_j$, we have
  \begin{equation*}
    d(\var\theta_{ij}-\deltacheck\eta_{ij})=0\,,
  \end{equation*}
  and we deduce, using Poincar\'e Lemma, that
  \begin{equation*}
    \var\theta_{ij}-\deltacheck\eta_{ij}=d\lambda_{ij}\,,
  \end{equation*}
  for $(\lambda_{ij})\in \Check{C}^1(\cover{U}_X, \sheaf{A}_X)$.
  An explicit calculation using \eqref{eq: theta new} and
  \eqref{eq:var_eta} confirms this relation with
  \begin{equation}
    \label{eq:var_lambda}
    \lambda_{ij} = 2\,\frac{w^{\prime\prime}_{ij}}{w^\prime_{ij}}
    \circ f_j\, \var f_j -\bigl( \log w^\prime_{ij}\circ f_j
    +\log z^\prime_{ij}\bigr)\, \var\log\del f_j -
    \Tilde\tau_{ij}\circ f_j\,\var f_j\,.
  \end{equation}
  The last term in this formula is obtained by varying the
  difference $f^*(\Tilde\tau_{ij})-\tau_{ij}$, that enters
  equation~\eqref{eq: theta new}. Clearly, the variation of
  $\tau_{ij}$ is zero and for the variation of
  $f^*(\Tilde\tau_{ij})$ we have
  \begin{align*}
    \var f^*(\Tilde\tau_{ij})
    &= \var \bigl( \Tilde\tau_{ij}\circ f_j\, df_j\bigr)
    = \var (\Tilde\tau_{ij}\circ f_j)\,df_j + 
    \Tilde\tau_{ij}\circ  f_j\, \var df_j\\
    &= {\Tilde\tau}^\prime_{ij}\circ f_j\,\var
    f_j\, df_j +\Tilde\tau_{ij}\circ f_j\,d\var f_j\\
    &= d\bigl( \Tilde\tau_{ij}\circ f_j\,\var f_j\bigr)\,,
  \end{align*}
  since $\Tilde\tau_{ij}\in\Omega^1(\Tilde U_i\cap\Tilde U_j)$ 
  (see section~\ref{sec:Theta}).

  Computing the coboundary of \eqref{eq:var_lambda} yields
  \begin{multline*}
    \deltacheck\lambda_{ijk} = -(\Tilde c_{ijk} + c_{ijk})
    \var\log\del f_k\\ -(\log w^\prime_{ij}\circ f_j +\log
    z^\prime_{ij})\,\var\log w^\prime_{jk}\circ f_k
    -(\deltacheck\Tilde\tau)_{ijk}\circ f_k\,\var f_k\,.
  \end{multline*}
  On the other hand, the variation of \eqref{eq:Theta} gives
  \begin{align*}
    \var\Theta_{ijk} &= \Tilde\phi^\prime_{ijk}\circ f_k\, \var
    f_k -(\Tilde c_{ijk} + c_{ijk})\,\var\log\del f_k -\log
    z^\prime_{ij}\, \var\log w^\prime_{jk}\circ f_k\\ &=
    \deltacheck\lambda_{ijk} +\log w^\prime_{ij}\circ
    f_j\;\var(\log w^\prime_{jk}\circ f_k)
    +\Tilde\phi^\prime_{ijk}\circ f_k\, \var f_k +(\deltacheck
    \Tilde\tau)_{ijk}\circ f_k\,\var f_k\,.
  \end{align*}
  Using the first equation in~\eqref{eq:tame_trivial}:
  \begin{equation*}
    d\Tilde\phi_{ijk} + (\deltacheck\tau)_{ijk} = -\log
    w^\prime_{ij}\, d\log w^\prime_{jk},
  \end{equation*}
  we get
  \begin{equation*}
    \var\Theta_{ijk}=\deltacheck\lambda_{ijk}\,.
  \end{equation*}
  Finally, putting it all together, we obtain
  \begin{equation*}
    \var\Omega = (a_i) +D(\eta + \lambda)\,,
  \end{equation*}
  as wanted.
\end{proof}
\begin{proof}[Second proof]
  The $2$-form $a_i$ in the relation $\var \omega_i = a_i
  +d\eta_i$ is a \emph{source} form~\cite{zuckerman}, hence it is
  uniquely determined by the de Rham class of $\omega_i$.
  Moreover, given a specific $\omega_i$, the form $d\eta_i$ is
  also determined (so $\eta_i$ is determined up to an exact
  form). Since $\omega_j = \omega_i + d\theta_{ij}$, we must have
  $a_i=a_j$ as both $a_i$ and $a_j$ are source forms for the
  \emph{same} Lagrangian problem. Here the requirement that the
  variation be vertical is crucial in order to ensure that $\var
  f/\del f$ glue as a geometric object --- a vector field on $X$.
  Therefore, from $\var\deltacheck\omega_{ij}
  =\deltacheck\var\omega_{ij}$, we get
  \begin{equation*}
    \var\theta_{ij}=\deltacheck\eta_{ij} + d\lambda_{ij}
  \end{equation*}
  by Poincar\'e lemma. Proceeding in the same fashion we also get
  \begin{equation*}
    d(\var\Theta_{ijk} - \deltacheck\lambda_{ijk})=0\,.
  \end{equation*}
  Now, both $\var\Theta_{ijk}$ and $\deltacheck\lambda_{ijk}$ are
  forms of degree one in the field direction, i.e. they contain
  one variation. Takens' acyclicity
  theorem~\cite{takens,zuckerman,del-freed} asserts the
  variational bicomplex is acyclic in all degrees except the top
  one in the de Rham direction, provided the degree in the
  variational direction is at least one. Hence,
  \begin{equation*}
    \var\Theta_{ijk} = \deltacheck\lambda_{ijk},
  \end{equation*}
  and we reach the same conclusion as in the previous proof.
\end{proof}

\subsection{Relative projective structures}
\label{sec:projective}

Here we interpret of the Euler-Lagrange equation from the
previous section through the principal $\gggg$-bundle over the
universal family of projective structures. First, we reformulate
theorem \ref{thm:variation} as follows
\begin{theorem}\label{thm:variation_hubbard}
  The Euler-Lagrange equation
  \begin{equation*}
    \delb\, h_i=\diff{D}_h\,\mu_i\,
  \end{equation*}
  for the vertical variational problem is the condition that
  push-forward of the projective connection $\{h_i\}$ onto
  $\Tilde X$ by the map $f$ is holomorphic.
\end{theorem}
\begin{proof} Indeed, the push-forward of $h$ is $f_*(h)=
\{\Tilde h\circ f_i^{-1}\cdot(\del f_i^{-1}/\del w_i)^2\}$, where
$\Tilde h_i = h_i - \{f_i,z_i\}$. It is a projective connection
on $\Tilde X$ because of the transformation law
\begin{equation*}
  \Tilde h_j - \Tilde h_i\circ z_{ij}(z^\prime_{ij})^2 =
  \{w_i,w_j\} \circ f_j \, (\del f_j)^2\,.
\end{equation*}
The Euler-Lagrange equation is equivalent to the equation
$\delb_{\mu}\Tilde h_i=0$, which is precisely the condition that
projective connection $f_*(h)$ is holomorphic on $\Tilde X$.
\end{proof}

It is well-known (see, e.g.~\cite{gunning}) that a holomorphic
projective connection on $X$ determines a projective structure on
$X$, and vice versa. The space of all projective structures on
$X$ is an affine space modeled over
$H^0(X,\sheaf{\Omega}^{\otimes 2}_{X})$ --- the vector space of
holomorphic quadratic differentials on $X$.

For any holomorphic family $C\to S$ of Riemann surfaces
parameterized by a complex manifold $S$, there is the holomorphic
family $P_S(C)\to S$ of relative projective structures on
$C$~\cite{hubbard}. The fiber over $s\in S$ is the affine space
of all (holomorphic) projective structures for $C_s$. We will be
interested in the universal case $S=\ttt$ and denote by $\ppp$
the universal family of relative projective structures.
Following~\cite{nag}, consider the following pullback diagram
\begin{equation}
\label{eq:proj_family}
\begin{CD}
   \sss @>>> \bbb\\ @VVV @VVV\\ \ppp @>>> \ttt
  \end{CD}
\end{equation}
where the vertical arrows are principal $\gggg$-bundles, and the
horizontal ones are affine bundles with spaces affine over
$H^0(X_\mu,\sheaf{\Omega}^{\otimes 2}_{X_\mu})$ as fibers,
$\mu\in\bbb$. (The curve $X_\mu$ depends only on the class of
$\mu$ modulo $\gggg$ and so do its holomorphic objects.) Here
$\sss$ is the space of all projective structures on $X$
holomorphic with respect to some complex structure determined by
$\mu\in\bbb$, without considering the quotient by $\gggg$. Since
every projective structure determines a complex structure, there
is an obvious projection $\sss \to \bbb$. As it follows from
theorem \ref{thm:variation_hubbard},
\begin{equation*}
  \sss = \{ (h,\mu)\in \mathcal{Q}(X)\times \bbb\vert
  \diff{D}_h\mu=\delb h\}\,,
\end{equation*}
so that $\sss$ is the critical manifolds for the mapping
$A:\mathcal{Q}(X)\times \bbb\to\CC^*$ (as well as for the map
$S:\mathcal{Q}(X)\times \bbb\to\CC/\ZZ(3)$).  The projection
$\sss \to \bbb$ is now just the projection on the second factor,
and every fiber over $\mu$ in $\sss$ is indeed an affine space
over the vector space $H^0(X_\mu,\sheaf{\Omega}^{\otimes
  2}_{X_\mu})$ (or rather its pull-back by $f$). Indeed, if
$(h,\mu)$ and $(h^\prime ,\mu)$ are two projective structures
subordinated to $\mu$, then we have
\begin{equation*}
 \diff{D}_h\mu=\delb h \qquad
 \text{and}\qquad\diff{D}_{h^\prime}\mu =\delb h^\prime\,,
\end{equation*}
which using the identity $\diff{D}_h\mu-\delb
h=\delb_{\mu}(\{f,z\}-h)$ imply that
\begin{equation*}
 \delb_\mu(h^\prime -h)=0\,,
\end{equation*}
concluding that $h-h^\prime$ is a $\mu$-holomorphic quadratic
differential.

The local meaning of the Euler-Lagrange equation --- the
condition $\diff{D}_h\mu=\delb h$ --- is the following.
\begin{lemma}
  The operators $\delb_\mu$ and $\diff{D}_h$ commute if and only
  if the Euler-Lagrange equation is satisfied.
\end{lemma}
\begin{proof}
  Let $v$ be a local section of $\sheaf{A}^{-1,0}_X$. As a result
  of a direct calculation we have (omitting the coordinate index
  $i$)
\begin{equation}
 \label{eq:commutation}
 \diff{D}_h\delb_\mu\,v - \delb_\mu\diff{D}_h\,v =
 L_v(\diff{D}_h\mu - \delb h)\,,
 \end{equation}
 where $L_v = v\,\del + 2\,\del v$ is the Lie derivative operator
 on $\sheaf{A}^{2,1}_X$. Thus the ``if'' part is clear.  For the
 ``only if'' part, assume the RHS of \eqref{eq:commutation} is
 zero for all $v$. Therefore, if we consider $f\,v$ for any local
 $f$, then we must have $v(f)\cdot (\diff{D}_h\mu - \delb h)=0$,
 implying \eqref{eq:proj_condition}.
\end{proof}

We conclude that $\sss$ is the geometric locus where the
commutativity condition $\diff{D}_h\delb_\mu=\delb_\mu\diff{D}_h$
is satisfied. Then we can consider $\diff{D}_h$ as a map between
two augmented Dolbeault complexes (where $\sheaf{\Theta}_\mu$ and
$\sheaf{\Omega}^{\otimes 2}_{X_\mu}$ are actually pull-backs of
the corresponding sheaves from $X_\mu$ to $X$ by the map
$f(\mu)$):
\begin{equation}
\label{eq:dolb}
\begin{CD}
  0@>>>\sheaf{\Theta}_\mu@>\imath>> \sheaf{A}^{-1,0}_X
  @>\delb_\mu>> \sheaf{A}^{-1,1}_X @>>> 0\\ & & @VV\diff{D}_hV
  @VV\diff{D}_hV @VV\diff{D}_hV \\ 0 @>>> \sheaf{\Omega}^{\otimes
    2}_{X_\mu} @>>\imath> \sheaf{A}^{2,0}_X @>>\delb_\mu>
  \sheaf{A}^{2,1}_X @>>> 0
\end{CD}
\end{equation}
where the morphism $\sheaf{\Theta}_\mu \xrightarrow{\diff{D}_h}
\sheaf{\Omega}^{\otimes 2}_{X_\mu}$ is now the usual third order
$\mu$-holomorphic operator~\cite{gunning-paper,hubbard}, also
familiar from the theory of the KdV equation~\cite{magri}. It
fits into the exact sequence
\begin{equation}
\label{eq:eich_sheaf}
\begin{CD}
  0 @>>> \sheaf{V}_X(h) @>\imath>> \sheaf{\Theta}_\mu
  @>\diff{D}_h>> \sheaf{\Omega}^{\otimes 2}_{X_\mu} @>>> 0
\end{CD}
\end{equation}
where $\sheaf{V}_X(h)$ is a rank three local system depending on
the projective structure $h$ --- a locally constant sheaf on $X$.
Actually, it is the sheaf of polynomial vector fields of degree
not greater than two in the coordinates adapted to $(h,\mu)$.
Passing to cohomology, we get:
\begin{equation}
\label{eq:eich_cohomology}
0 \to H^0(X_\mu,\sheaf{\Omega}^{\otimes 2}_{X_\mu}) \to
H^1(X,\sheaf{V}_X(h)) \to H^1(X_\mu,\Theta_\mu) \to 0
\end{equation}
According to the theorem of Hubbard~\cite{hubbard}, sequence
\eqref{eq:eich_cohomology} is isomorphic to the tangent bundle
sequence for the relative projective structure $\ppp \to \ttt$ at
$(h,\mu)$.  Furthermore, the usual machinery of local systems
shows that $H^1(X,\sheaf{V}_X(h))$ is isomorphic to the Eichler
cohomology group $H^1(\pi_1(X,p),\boldsymbol{V}(h)_p)$, where
$\boldsymbol{V}(h)_p$ is the stalk of $\sheaf{V}_X(h)$ over the
point $p$. The proof that this coincides with the classical
Eichler cohomology (see~\cite{kra}), can be obtained by lifting
everything to the universal cover $\HH$ of $X$ and using factors
of automorphy (see~\cite{hubbard} for further details).

On the other hand, from our description of $\sss$ we have
\begin{equation}
  \label{eq:proj_condition_diff}
  T_{(h,\mu)}\sss =
  \{ (\dot h,\dot\mu)\in A^{2,0}(X)\times
  A^{-1,1}(X)\vert \diff{D}_h\dot\mu=\delb_\mu\dot h\}
\end{equation}
and the RHS can be written as the fiber product
\begin{equation*}
  A^{2,0}(X)\times_{A^{2,1}(X)}A^{-1,1}(X)
\end{equation*}
with respect to the pair of maps $\delb_\mu$ and $\diff{D}_h$.
For vertical --- along the fiber of $\sss\to \bbb$ --- tangent
vectors to $\sss$ at $(h,\mu)$ we have
\begin{equation*}
  (\dot h,\dot \mu) = (\diff{D}_h v,\delb_\mu v) \,,
\end{equation*}
where $v\in A^{-1,0}(X)$ is the infinitesimal generator. This
pair clearly satisfies the condition in
\eqref{eq:proj_condition_diff}, since $\sss$ is the geometric
locus of the commutativity condition.  Thus the map sending $v
\mapsto(\diff{D}_h v,\delb_\mu v)$ describes the vertical tangent
bundle of $\sss\to \bbb$. Therefore, if $[h,\mu]$ denotes the
class of $(h,\mu)$, we have for the vertical tangent space to
$\sss$ at $[h,\mu]$:
\begin{equation}\label{eq:TP}
  T_{V,[h,\mu]}\sss\cong
  \bigl(A^{2,0}(X)\times_{A^{2,1}(X)}A^{-1,1}(X)\bigr) /
  (\diff{D}_h,\delb_\mu )(A^{-1,0}(X))\,,
\end{equation}
which obviously projects onto $H^{-1,1}_{\delb_\mu}(X)\cong
H^1(X_\mu,\Theta_\mu)$. Now, this is just the $C^\infty$ image of
the Eichler cohomology description of the tangent sheaf to the
relative projective structure $\ppp\to \ttt$ and we have the
following
\begin{proposition}
The differential geometric description of the tangent space to
$\ttt$ at the class of $(h,\mu)$ as given by \eqref{eq:TP}
coincides with the algebraic description given by the Eichler
cohomology group $H^1(X,\sheaf{V}_X(h))$.
\end{proposition}
\begin{proof}
  Consider the cone of $\diff{D}_h:\sheaf{A}^{-1,\bullet}_X \to
  \sheaf{A}^{2,\bullet}_X$:
  \begin{equation*}
    \sheaf{C}^\bullet_X: 0 \longrightarrow \sheaf{A}^{-1,0}_X
    \xrightarrow{\delb_\mu\oplus\diff{D}_h}
    \sheaf{A}^{-1,1}_X\oplus
    \sheaf{A}^{2,0}_X \xrightarrow{\diff{D}_h -\delb_\mu}
    \sheaf{A}^{2,1}_X \longrightarrow 0
  \end{equation*}
  Its cohomology sheaf complex equals $\sheaf{V}_X(h)$, thus by
  standard homological algebra arguments (see,
  e.g.~\cite{maclane}) one has
  $\HHH^1(X,\sheaf{C}^\bullet_X)=H^1(X,\sheaf{V}_X(h))$ and from
  the canonical sequence
  \begin{equation*}
    0 \longrightarrow \sheaf{A}^{2,\bullet}_X[-1] \longrightarrow
    \sheaf{C}^\bullet_X \longrightarrow \sheaf{A}^{-1,\bullet}_X
    \longrightarrow 0
  \end{equation*}
  one gets \eqref{eq:eich_cohomology}. On the other hand, the RHS
  of \eqref{eq:TP} is the first cohomology group of the complex
  \begin{equation*}
    0 \longrightarrow A^{-1,0}(X)
    \xrightarrow{\delb_\mu\oplus\diff{D}_h} A^{-1,1}(X)\oplus
    A^{2,0}(X) \xrightarrow{\diff{D}_h -\delb_\mu}
    A^{2,1}(X)\longrightarrow 0
  \end{equation*}
  which is equal to the first term
  \begin{equation*}
    {}^\backprime E^{p,q}_1 \cong \Check{H}^q(X,\sheaf{C}^p_X) =
    \begin{cases}
      C^p(X) & q=0\,,\\ 0 & q>0\,.
    \end{cases}
  \end{equation*}
  of the spectral sequence computing
  $\HHH^\bullet(X,\sheaf{C}^\bullet_X)$.
\end{proof}

\subsection{Geometry of the vertical variation}
\label{sec:vert_variation}

Here we consider functional $A[f]$ as as map
$A:\mathcal{Q}(X)\times \bbb\longrightarrow \CC^*$, where
$\mathcal{Q}(X)$ is the affine space of all $C^\infty$ projective
connections on $X$ and $\bbb$ is the total space of the
Earle-Eells fibration.

By theorem~\ref{thm:variation_hubbard}, the critical manifold for
$A[f]$ coincides with $\sss$. Considering critical values of $A$
(''on shell" condition) leads to the function $A:
\sss\longrightarrow\CC^*$, where $A(h,\mu)=\langle
\Omega[h,\mu],\Sigma\rangle_m$. Since $\sss$ is a principal
$\gggg$-bundle over $\ppp$, it is interesting to analyze the
behavior of $A$ under the $\gggg$-action.  It is given by the
following
\begin{lemma}
  The directional derivative of the action functional $A$ for the
  vertical tangent vector $(\diff{D}_h v,\delb_\mu v)$ to
  $\sss$ at $(h,\mu)$, where $v\in A^{-1,0}(X)$, is given by
  \begin{equation} \label{eq:g_diff}
    4\pi\sqrt{-1} \int_X \mu\,\diff{D}_h v \cdot A.
  \end{equation}
\end{lemma}
\begin{proof}
  We just repeat the computation of the vertical variation with
  additional term $2\,\mu\,\var h\, dz\wedge d\Bar z$, where
  $\var h=\diff{D}_h v$.  Since the main term, given by
  $2\,(\delb h -\diff{D}_h\mu)\,v \,dz\wedge d\Bar z$ vanishes
  "on shell", this proves the result.
\end{proof}
Formula \eqref{eq:g_diff} defines a function $c: \sss
\longrightarrow \bigl(\mathrm{Lie}\;\gggg\bigr)^*$ by assigning
to the pair $(h,\mu)$ a linear functional on
$\mathrm{Lie}\;\gggg\cong A^{-1,0}(X)$ as follows:
\begin{equation}
v\mapsto 2\int_X \mu\,\diff{D}_h v\,.
\end{equation}
Equivalently, $c$ is a 1-cochain over $\mathrm{Lie}\;\gggg$ with
values in functions over $\sss$ with left
$\mathrm{Lie}\;\gggg$-action.
\begin{proposition}\label{prop:cocycle}
  The 1-cochain $c$ is a 1-cocycle.
\end{proposition}
\begin{proof}
  For $v,w\in A^{-1,0}(X)\cong \mathrm{Lie}\;\gggg$ we have
  \begin{equation*}
    \delta c(v,w) = v\cdot c(w) - w\cdot c(v) - c([v,w])
  \end{equation*}
  where $c(u):\sss\to \CC$ is the function
  \begin{equation*}
    c(u)(h,\mu)=2\int_X \mu\,\diff{D}_h u\,.
  \end{equation*}
  Using the infinitesimal action,
  \begin{gather*}
    v\cdot c(w)(h,\mu) = 2\int_X \bigl(%
    \delb_\mu\,v\,\diff{D}_h w +
    \mu\,L_v(\diff{D}_h w\bigr)\,,\\ \intertext{we get}
    \begin{split}
    (\delta c)(h,\mu) = 2\int_X \bigl(%
    \delb_\mu\,v\,\diff{D}_h w &+ \mu\,L_v(\diff{D}_h w \\
    -\delb_\mu\,w\,\diff{D}_h v &-\mu\,L_w(\diff{D}_h v)
    -\mu\,\diff{D}_h\,L_v w \bigr)\,.
    \end{split}
  \end{gather*}
  where $L_v=v\del+2\del v$ is the Lie derivative on $A^2(X)$,
  and the Lie bracket in $A^{-1,0}(X)$ is the usual vector field
  Lie bracket: $[v,w] = L_vw = (v\,\del w - w\,\del v)$. Using
  the identity
  \begin{equation*}
    L_v(\diff{D}_h w) - L_w(\diff{D}_h v) - \diff{D}_h L_v (w)=0,
  \end{equation*}
  we are left with
  \begin{align*}
    (\delta c)(v,w)(h,\mu) &= 2\int_X \bigl(%
    \delb_\mu v\,\diff{D}_h w -\delb_\mu w\,\diff{D}_h v
    \bigr)\\
    &= 2\int_X v\,\bigl(%
    \diff{D}_h\delb_\mu w -\delb_\mu\diff{D}_h w \bigr)\\ &=
    0\,,
  \end{align*}
  because of the commutativity condition and the skew-symmetry of
  the operator $\diff{D}_h$.
\end{proof}

\appendix

\section{Appendix}
\label{sec:appendix}

\subsection{Cones}
\label{sec:cones}
Recall~\cite{maclane} that for a map $u:\AAA^\bullet
\to\BBB^\bullet$ the \emph{cone} $\CCC^\bullet_u$ of $u$ is the
complex:
\begin{equation*}
  \CCC^\bullet_u = \AAA^\bullet[1]\oplus \BBB^\bullet
\end{equation*}
with differential
\begin{equation*}
  d(a,b)=(-da,u(a) + db)\,.
\end{equation*}
The cone fits into the exact sequence:
\begin{equation*}
  0\longrightarrow \BBB^\bullet \longrightarrow \CCC^\bullet_u
  \longrightarrow \AAA^\bullet[1] \longrightarrow 0\,.
\end{equation*}
If the map $u$ is injective, this is the same as the cokernel of
$u$ (up to a shift in the resulting exact cohomology
sequence). 

For the Deligne complex, we often find that the equivalent
definition~\cite{esnault,del-freed} of
$\ZZ(p)^\bullet_\mathcal{D}$ is
\begin{equation}
  \label{eq:del_cone}
  \ZZ(p)^\bullet_\mathcal{D}=\mathrm{Cone}\bigl( \ZZ(p)\oplus
  F^p(\sheaf{A})^\bullet_M \xrightarrow{\imath - \jmath}
  \sheaf{A}^\bullet_M\bigr)[-1]\,,
\end{equation}
where $\jmath :F^p(\sheaf{A})^\bullet_M \to \sheaf{A}^\bullet_M$
is the Hodge-Deligne filtration (\emph{filtration b\^ete}), that
is, the $n$-th sheaf of $F^p(\sheaf{A})^\bullet_M$ is
$\sheaf{A}^n_M$ if $n\geq p$, and zero otherwise.

Briefly, the equivalence is shown as follows. The cone in
\eqref{eq:del_cone} is equal to
\begin{equation*}
  \mathrm{Cone}\bigl(\ZZ(p)\longrightarrow
  \mathrm{Cone}(F^p(\sheaf{A})^\bullet_M \longrightarrow
  \sheaf{A}^\bullet_M) \bigr)[-1]\,.
\end{equation*}
The inner cone can clearly be replaced by the cokernel of the
inclusion map, namely the (sharp) truncation $\tau^{\leq
  p-1}\sheaf{A}^\bullet_M$ of the de Rham complex. Thus we have
\begin{equation*}
  \mathrm{Cone}\bigl( \ZZ(p)\longrightarrow \tau^{\leq
  p-1}\sheaf{A}^\bullet_M \bigr)[-1]\,,
\end{equation*}
which equals $\ZZ(p)^\bullet_\mathcal{D}$ as defined in the main
text.

\subsection{Fundamental class}
\label{sec:app_fund_class}
We want to collect here some technical facts and computations
related to the construction of a representative of the
fundamental class $[M]$, that are not strictly necessary in the
main body of this paper.

Recall that we work with the double complex
\begin{equation*}
  \SSS_{p,q}=S_p(N_q(U\to X))\,,
\end{equation*}
where $N_\bullet(U\to X)$ is the nerve of the covering $U\to X$,
and $S_\bullet$ is the singular simplices functor.

\subsubsection{}
We saw in the main text, sec.~\ref{sec:prel_fundamental}, that
when $\SSS_{p,\bullet}$ \emph{resolves} $S_p(M)$ for any fixed
$p$, the total homology of $\SSS_{\bullet,\bullet}$ is equal to
$H_\bullet(M,\ZZ)$. By definition, this condition is that
$H_0(S_p(N_\bullet(U))) \cong S_p(M)$ and
$H_q(S_p(N_\bullet(U)))=0$ for $q>0$. Then the isomorphism
$H_\bullet(M,\ZZ) \cong H_\bullet(\Tot \SSS)$ can be easily
obtained by carefully lifting a cocycle in $S_\bullet(M)$ to a
total cocycle in $\SSS_{\bullet,\bullet}$.\footnote{See, e.g.,
\protect\cite{maclane}. Details for this calculation can be found
in the appendix of \protect\cite{prev_paper}.}  More concisely,
we have ${}^{\backprime}E^1_{p,q}=H_q(S_p(N_\bullet(U)))=0$ for
$q>0$ (the spectral sequence collapses) and at the next step one
has ${}^{\backprime}E^2_{p,0} = {}^{\backprime}E^\infty_{p,0}
\cong H_p(S_\bullet(M)) = H_p(M,\ZZ)$, as wanted.

These requirements are met for a \v Cech covering $\cover{U}_M$,
where a contracting homotopy for $S_p(N_\bullet(\cover{U}_M))$
can be constructed explicitly~\cite{weil} (see
also~\cite{goldberg}, appendix on the de Rham theorem). Indeed,
one can easily show that $H_0(S_p(N_\bullet(\cover{U}_M))\cong
S_p(M)$ by applying $S_p(-)$ to the sequence $\cdots
N_1(\cover{U}_M)\rightrightarrows N_0(\cover{U}_M)\rightarrow M$.
The resulting maps are $\coprod_i\sigma_i \to \sum_i\sigma_i$ and
$\coprod_{ij}\sigma_{ij}\to \coprod_i(\sum_j (\sigma_{ji} -
\sigma_{ij}))$, so the composition is zero. Moreover, if
$\sum_i\sigma_i = 0$, for any pair of indices $ij$, we must have
$\sigma_i\rvert_{U_{ij}} + \sigma_j\rvert_{U_{ij}}=0$, so that
$\sigma_i\rvert_{U_{ij}}\coprod
\sigma_j\rvert_{U_{ij}}=\sigma_i\rvert_{U_{ij}} \coprod -
\sigma_i\rvert_{U_{ij}}= \delpp \sigma_i\rvert_{U_{ij}}$, proving
the claim.

Similarly, if $U\to M$ is a regular covering with
$G=\mathrm{Deck}(U/M)$ acting on the right on $U$, then
$\SSS_{p,0}=S_p(N_0(U))\equiv S_p(U)$ is a free (right)
$G$-module~\cite{maclane}, so that $S_p(N_\bullet(U)) \cong
S_p(U)\otimes_{\ZZ G}B_\bullet(G)$ resolves $S_p(U)\otimes_{\ZZ
G}\ZZ \cong S_p(M)$ hence
\begin{equation*}
  H_q(S_p(N_\bullet(U\to M)))\cong
  \begin{cases}
    S_p(M)& q=0\\ 0 & q >0\,,
  \end{cases}
\end{equation*}
as wanted.

\subsubsection{}
\label{sec:abstract_nerve}
Since $\SSS_{\bullet,\bullet}$ is a double complex, it is well
known that its associated total complex can be filtered in two
ways --- with respect to either $p$ or $q$.  Filtering over the
second index of $\SSS_{p,q}=S_p(N_q(U))$ yields the second
spectral sequence with
\begin{equation*}
  {}^{\backprime\backprime}E^1_{p,q} \cong H^{\delp}_p
  (\SSS_{\bullet,q}) \equiv H_p(S_\bullet(N_q(U))\,.
\end{equation*}
Although not required in the following it is interesting to see
when and whether this latter sequence also degenerates, like the
other one.  In other words, we want to consider the case when for
fixed $q$ the complex $\SSS_{\bullet,q}$ is acyclic in degree
$>0$.

\begin{assumption}
  The covering $U\to M$ is \emph{good}, that is, each
  $N_q(U)=U\times_M\dots\times_M U$ is contractible, hence is
  acyclic for the singular simplices functor.
\end{assumption}
\begin{remark}
  The assumption on $U\to M$ guarantees the de Rham complex is a
  resolution of $\CC$, so the second cohomological spectral
  sequence $H^p(\Check{C}^q(U;\sheaf{A}^\bullet))$ degenerates
  and the total cohomology equals $\Check{H}^q(U;\CC)$.
\end{remark}
By virtue of the assumption, ${}^{\backprime\backprime}E^1$ is
computed as
\begin{equation*}
  {}^{\backprime\backprime}E^1_{q,p} \cong
  \begin{cases}
    \ZZ <N_q \cat{R}_U> & p=0\\ 0 & p>0\,,
  \end{cases}
\end{equation*}
where $N_q \cat{R}_U$ is the set of connected components of
$N_q(U)$ and $\ZZ <N_q \cat{R}_U>$ is the abelian group generated
by $N_q \cat{R}_U$.  This follows from the fact that $H_0$ gives
us a factor $\ZZ$ for every connected component of
$N_q(U)$. These connected components arrange into a simplicial
set $N_\bullet \cat{R}_U$, where the face maps are induced by the
face maps of the nerve $N_\bullet(U)$, specifying where every
component goes. Thus $N_\bullet \cat{R}_U$ expresses the pure
combinatorics of the covering.  Since the spectral sequence
collapses, the total homology is equal to
\begin{gather*}
  {}^{\backprime\backprime}E^2_{q,0}=
  {}^{\backprime\backprime}E^\infty_{q,0} \cong H_q(\ZZ<N_\bullet
  R_U>)\\ \intertext{and (see~\cite{may})} H_q(\ZZ <N_\bullet
  R_U>)\cong H_q(\lvert N_\bullet R_U\rvert)\,,
\end{gather*}
where $\lvert\cdot\rvert$ is the geometric realization of
$N_q\cat{R}_U$, namely, the CW-complex obtained by putting in a
standard $q$-simplex $\Delta^q$ for each element in
$N_q\cat{R}_U$ and gluing them together according to the face
maps. Therefore, for a good covering the three homologies are
equal:
\begin{equation*}
  H_q(\Tot \SSS_{\bullet,\bullet})\cong H_q(M,\ZZ) \cong
  H_q(\lvert N_\bullet \cat{R}_U\rvert)\,.
\end{equation*}
In our concrete examples, an ordinary \v{C}ech covering is good
if all $U_{i_0}\cap \dots \cap U_{i_q}$ are contractible. In this
case, to compute $H^{\delp}_0(S_\bullet (N_q\cover{U}_M))$ we
must assign a $\ZZ$ factor to each $U_{i_0}\cap \dots \cap
U_{i_q}$.  Following~\cite{weil,goldberg}, denote $U_{i_0}\cap
\dots \cap U_{i_q}$ as a generator in this group by the symbol
$\Delta_{i_0,\dots,i_q}$, so that $\ZZ<N_q \cat{R}_U> =
\bigoplus_{i_0,\dots ,i_q}\ZZ \cdot \Delta_{i_0,\dots,i_q}$ and
$N_q \cat{R}_U = \{\Delta_{i_0,\dots,i_q}\}$. Therefore,
$N_\bullet \cat{R}_U$ represents the abstract nerve of the open
cover and $\lvert N_\bullet \cat{R}_U\rvert$ is the CW-complex
obtained by replacing each $\Delta_{i_0,\dots,i_q}$ --- in other
words, each non void intersection --- by a standard $q$-simplex
and gluing them according to the face maps of $N_\bullet
\cat{R}_U$.

On the other hand, if $U\to M$ is a $G$-covering, then according
to~\ref{sec:prel_coverings} $N_q(U)= U\times G^q$, and it is good
if $U$ is contractible. Thus $\ZZ <N_\bullet \cat{R}_U> \cong
\ZZ\otimes_{\ZZ G}B_\bullet (G)$, so that
\begin{equation*}
  H_q(\ZZ <N_\bullet \cat{R}_U>) \cong H_q(G;\ZZ) \cong
  H_q(BG,\ZZ)\,,
\end{equation*}
where $BG=\lvert N_\bullet \cat{R}_U\rvert$ is the classifying
space of $G$, where in this case $N_q \cat{R}_U = G^q$ for $q\geq
1$ and $N_0\cat{R}_U = \mathrm{point}$.

\subsubsection{}
\label{sec:fund_class_computation}
Let us return to the main problem of representing the fundamental
class of $X$ as a total cycle in the double complex $\SSS_{p,q}$.
If the sequence $0\leftarrow S_p(X)\leftarrow S_p(N_\bullet (U))$
is exact, then there exists a splitting $\tau:S_p(X) \to
S_p(N_\bullet (U))$, i.e.~the map $\tau$ satisfies $\epsilon\circ
\tau=\id_{S_p(X)}$. In other words, $\tau$ is the first step of
an explicit contracting homotopy for $S_p(N_\bullet (U))$.  Then
a cycle representing $[X]$ can be produced by lifting $X$ via
$\tau$ and completing $\tau(X)$ to a total cycle using the
standard descent argument.

In the concrete examples we have been looking at, this can be
done as follows. The case where $U$ is a regular $G$-covering can
be handled by starting from a fundamental domain $F$ for the
action of $G$ on $U$, where we regard $F$ as an element of degree
$(p,0)$ in $\SSS_{p,0}\cong S_p(U)\otimes_{\ZZ G}B_0(G)\cong
S_p(U)$. Full details are spelled out in~\cite{prev_paper}. If
$U$ comes from an ordinary \v Cech covering $\cover{U}_X$, we
first replace $S_p(X)$ by $\cover{U}_X$-small simplices:
\begin{equation*}
  0\longleftarrow S^\cover{U}_p(X) \longleftarrow S_p(N_\bullet
  \cover{U}_X),
\end{equation*}
where the $\cover{U}_X$-small simplices are those whose support
is contained in the open cover $\cover{U}_X=\{U_i\}$. Second,
write $X=\sum_i \sigma_i$, where all $\sigma_i$ are
$\cover{U}_X$-small, and set $\tau (X)=\sum_i \sigma_i\cdot
\Delta_i \eqdef \Sigma_0\in \SSS_{p,0}$. Since
\begin{equation*}
  \epsilon (\delp\Sigma_0)=\delp\epsilon\Sigma_0
  =\delp\epsilon\tau (X) =\delp X \equiv 0\,,
\end{equation*}
by the standard argument there exist
$\Sigma_1,\Sigma_2,\dots,\Sigma_p$, with $\Sigma_k\in
\SSS_{p-k,k},\,k=1,\ldots, p,$ such that
\begin{equation*}
  \delp \Sigma_0 = \delpp \Sigma_1\,,\dots\,, \delp \Sigma_{q-1}
  = \delpp \Sigma_q\,,\dots\,, \delp \Sigma_p=0\,.
\end{equation*}
This schema can be implemented in a fairly explicit way using a
map $h:N_\bullet \cat{R}_{\cover{U}}\to S^\cover{U}_\bullet(X)$
constructed in~\cite{bott-tu} (Th. 13.4, proof) to realize the
nerve of a covering. Of course, our case of interest here is
$p=2$.

\begin{figure}[t]
  \begin{center}
    \includegraphics[scale=.3]{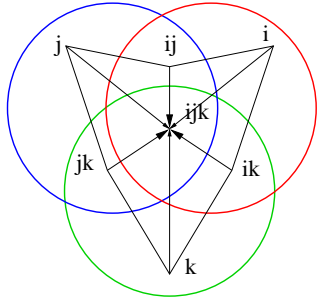}
  \end{center}
  \caption{\label{fig}Intersections and their nerve.}
\end{figure}

In order to describe $h$ we shall need the barycentric
decomposition $N_\bullet\Tilde{\cat{R}}_\cover{U}$ of
$N_\bullet\cat{R}_\cover{U}$ (see~\cite{segal} for a more
complete explanation). For any finite subset $\tau$ of the index
set $I$ denote $U_\tau= \cap_{i\in \tau}U_i$, and let:
\begin{gather*}
  N_0\Tilde{\cat{R}}_\cover{U} = \coprod_{\tau\subset
    I}\{U_\tau\}\\ N_1\Tilde{\cat{R}}_\cover{U} =
    \coprod_{\tau_0\subset \tau_1 \subset I} \{U_{\tau_1}\subset
    U_{\tau_0}\}\\ \dotsm\\ N_q\Tilde{\cat{R}}_\cover{U} =
    \coprod_{\tau_0\subset\dots\subset \tau_q\subset I}
    \{U_{\tau_q}\subset \dots\subset U_{\tau_0}\}\,.
\end{gather*}
In order to construct the mapping $h$, assign to each $U_\tau$ a
point $v_\tau\in U_\tau$, to any inclusion $U_{\tau_1}\subset
U_{\tau_0}$ a path from $v_{\tau_0}$ to $v_{\tau_1}$, and to
$U_{\tau_2}\subset U_{\tau_1}\subset U_{\tau_0}$ the cone from
$v_{\tau_0}$ to the path from $v_{\tau_1}$ to $v_{\tau_2}$, which
is of course a 2-simplex.  Denote by $\Delta (v_{\tau_0})$,
$\Delta (v_{\tau_0},v_{\tau_1})$ and $\Delta
(v_{\tau_0},v_{\tau_1},v_{\tau_2})$ the $0$, $1$ and
$2$-simplices so obtained. Observe how the simplices constructed
in this way inherit an orientation from the natural one on the
barycentric decomposition $N_\bullet\Tilde{\cat{R}}_\cover{U}$;
this is the main reason for using
$N_\bullet\Tilde{\cat{R}}_\cover{U}$ in place of
$N_\bullet{\cat{R}}_\cover{U}$. So, for example, $\Delta
(v_{\tau_0},v_{\tau_1},v_{\tau_2})$ has the orientation induced
by the order $v_{\tau_0}\leq v_{\tau_1}\leq v_{\tau_2}$
associated to the inclusion $\tau_0\subset \tau_1 \subset
\tau_2$.  The typical situation for the indices $i,j,k$ looks as
in figure~\ref{fig}: to the index sets $i$, $ij$ and $ijk$
correspond the points $v_i$, $v_{ij}$ and $v_{ijk}$ in $U_i$,
$U_{ij}$ and $U_{ijk}$, respectively.  Then $\Delta(v_i,v_{ij})$
is the $1$-simplex joining $v_i$ and $v_{ij}$,
$\Delta(v_{ij},v_{ijk})$ the one joining $v_{ij}$ and $v_{ijk}$,
and so on.
After these preparations, define an element $\Sigma_0$ in
$\SSS_{2,0}$ as
\begin{equation*}
  \Sigma_0 = \sum_{i\in I}\mathrm{st}(v_i)\cdot \Delta_i\,,
\end{equation*}
where
\begin{equation*}
  \mathrm{st}(v_i)=\sum_{j,k: \Delta_{ijk}\neq 0} \epsilon_{ijk}
  \bigl( \Delta (v_i,v_{ij},v_{ijk}) - \Delta
  (v_i,v_{ik},v_{ijk}) \bigr)
\end{equation*}
is the star of the vertex $v_i$, and $\epsilon_{ijk}=\pm 1$
according to whether the order of the triple $i,j,k$ agrees the
orientation or not, namely whether the order $ijk$ is the same as
the cyclic (counterclockwise) order around the vertex
$v_{ijk}$. Recall that $\Delta_\tau$ is the symbol corresponding
to $U_\tau$, when considered as a generator in the abelian group
generated by the nerve, as in~\ref{sec:abstract_nerve}. Rewriting
$\Sigma_0$ as
\begin{align*}
  \sum_{i\in I}\mathrm{st}(v_i)\cdot \Delta_i
  = \sum_{\langle i,j,k \rangle} \epsilon_{ijk} \biggl\{%
  &\bigl( \Delta (v_i,v_{ij},v_{ijk}) -
  \Delta(v_i,v_{ik},v_{ijk}) \bigr) \cdot \Delta_i \\ -&\bigl(
  \Delta (v_j,v_{ij},v_{ijk}) - \Delta (v_j,v_{jk},v_{ijk})
  \bigr) \cdot \Delta_j\\ +&\bigl( \Delta (v_k,v_{ik},v_{ijk}) -
  \Delta (v_k,v_{jk},v_{ijk}) \bigr) \cdot \Delta_k \biggr\},
\end{align*}
where $\sum_{\langle i,j,k\rangle}$ means sum over triples of
indices in $I$, its first differential is:
\begin{align*}
  \delp\Sigma_0 &= \sum_{\langle i,j,k\rangle} \epsilon_{ijk}
  \bigl\{ \Delta (v_{ik},v_{ijk})\cdot (\Delta_k -\Delta_i)\\
  &\qquad - \Delta (v_{ij},v_{ijk})\cdot (\Delta_j -\Delta_i)
  -\Delta (v_{jk},v_{ijk})\cdot (\Delta_k -\Delta_j) \bigr\} \\
  &+ \sum_{\langle i,j,k\rangle} \epsilon_{ijk} \bigl\{ \Delta
  (v_i,v_{ij})\cdot \Delta_i -\Delta (v_i,v_{ik})\cdot \Delta_i
  -\Delta (v_j,v_{ij})\cdot \Delta_j\\ &\qquad +\Delta
  (v_j,v_{jk})\cdot \Delta_j +\Delta (v_k,v_{ik})\cdot \Delta_k
  -\Delta (v_k,v_{jk})\cdot \Delta_k \bigr\}.
\end{align*}
The last sum is easily seen to be zero, while the first can be
rewritten as $\delpp\Sigma_1$ for the following element in
$\SSS_{1,1}$:
\begin{equation*}
  \Sigma_1 = \sum_{\langle i,j,k\rangle} \epsilon_{ijk} \bigl\{
  \Delta (v_{ik},v_{ijk})\cdot \Delta_{ik} -\Delta
  (v_{ij},v_{ijk})\cdot \Delta_{ij} -\Delta (v_{jk},v_{ijk})\cdot
  \Delta_{jk} \bigr\}
\end{equation*}
Again, computing the first differential gives
\begin{align*}
  \delp\Sigma_1 &= \sum_{\langle i,j,k\rangle} \epsilon_{ijk}
  \bigl\{ v_{ijk}\cdot (\Delta_{ik}-\Delta_{ij}-\Delta_{jk})
  \bigr\}\\ &+ \sum_{\langle i,j,k\rangle} \epsilon_{ijk} \bigl\{
  v_{ij}\cdot \Delta_{ij} -v_{ik}\cdot \Delta_{ik} +v_{jk}\cdot
  \Delta_{jk} \bigr\},
\end{align*}
with the last sum being identically zero. The first term can be
rewritten as $\delpp\Sigma_2$, where
\begin{equation*}
  \Sigma_2 = - \sum_{\langle i,j,k\rangle} \epsilon_{ijk}
  v_{ijk}\cdot \Delta_{ijk}\,.
\end{equation*}
Finally, the total chain $\Sigma\equiv \Sigma_0 + \Sigma_1
-\Sigma_2$ is a cycle, $\del \Sigma=0$, and we have the following
expression for the representative of the fundamental class of $X$
in the double complex:
\begin{align*}
  \Sigma &= \sum_{\langle i,j,k
    \rangle} \epsilon_{ijk} \biggl\{%
  \bigl( \Delta (v_i,v_{ij},v_{ijk}) - \Delta(v_i,v_{ik},v_{ijk})
  \bigr) \cdot \Delta_i \\ &\qquad\qquad\qquad -\bigl( \Delta
  (v_j,v_{ij},v_{ijk}) - \Delta (v_j,v_{jk},v_{ijk}) \bigr) \cdot
  \Delta_j\\ &\qquad\qquad\qquad +\bigl( \Delta
  (v_k,v_{ik},v_{ijk}) - \Delta (v_k,v_{jk},v_{ijk}) \bigr) \cdot
  \Delta_k \biggr\}\\ &+ \sum_{\langle i,j,k\rangle}
  \epsilon_{ijk} \bigl\{ \Delta (v_{ik},v_{ijk})\cdot \Delta_{ik}
  -\Delta (v_{ij},v_{ijk})\cdot \Delta_{ij} -\Delta
  (v_{jk},v_{ijk})\cdot \Delta_{jk} \bigr\}\\ &+ \sum_{\langle
  i,j,k\rangle} \epsilon_{ijk} v_{ijk}\cdot \Delta_{ijk}\,.
\end{align*}
\begin{remark}
  By taking the second augmentation, the total cycle $\Sigma$
  maps to:
  \begin{equation*}
    \sum_{\langle i,j,k\rangle} \epsilon_{ijk} \Delta_{ijk}\,,
  \end{equation*}
  which is the 2-cycle in the CW complex representing the
  combinatorics of the cover $\cover{U}$, and therefore the
  homology of $X$, in degree $p=2$.
\end{remark}

\section{Acknowledgements}
At the early stage of this work we appreciated useful discussions
with J.L.~Dupont and especially C.-H.~Sah, who passed away in
July 1997. His generosity of mind and enthusiasm made all our
discussions special. He is deeply missed.
  
The work of L.T. was partially supported by the NSF grant
DMS-98-02574.

\end{document}